# An approach to Set-Valued Anti-Homomorphism of $\mathcal{T}$-Rough Ideals on $\Gamma$-Semigroups

M.Thangeswari[1], R. Muthucumaraswamy[2], K. Anitha[3*]

[1] Research Scholar, Department of Mathematics,
Sri Venkateswara College of Engineering, Anna University, Chennai, India
[2] Department of Mathematics, Sri Venkateswara College of Engineering, Sriperumbudur, Chennai, India
[3*] Department of Mathematics, SRM Institute of Science and Technology, Ramapuram, Chennai, India.
Corresponding Author: anithak1@srmist.edu.in

**Abstract**

In this paper, the concepts of set-valued anti-homomorphism and strong set-valued anti-homomorphism of $\Gamma$-semigroup are introduced. The notions of generalized lower and upper approximation operators, constructed by means of set-valued mapping, which is a generalization of the notion of lower and upper approximation of $\Gamma$-semigroup, are provided. Some significant properties of lower and upper approximations of $\Gamma$-semigroup under (strong)set-valued anti-homomorphisms are discussed. We provide examples for the properties of lower and upper approximations of $\Gamma$-semigroup under (strong)set-valued anti-homomorphisms. This article explores the notion of a set-valued anti-homomorphism for $\Gamma$-semigroups, which is a generalization of anti-homomorphism in general.

## 1    Introduction

A mathematical framework for examining ambiguity is termed as Rough set [1-4]. Knowledge has recently received wide attention on the research areas in both of the real-life applications and the theory itself. It has found practical applications in many areas such as knowledge discovery, machine learning, data analysis, pattern recognition, approximate classification, conflict analysis, and so on. Rough set theory is a mathematical framework for dealing with uncertainty and to some extent overlapping fuzzy set theory. The rough set theory approach is based on indiscernibility relations and approximations. The theory of rough set is an extension of set theory. It is a revolutionary mathematical approach to dealing with imperfect data. By using a set's boundary area instead of its membership, rough set theory expresses contradiction.

Rough set theory includes three basic elements: the universe set, the binary relations and a subset described by a pair of ordinary sets. In the past few years, most studies have been focusing on the binary relations and the subsets; many interesting and constructive extensions to binary relations and the subsets. But few researchers have paid little attention to another basic element: the universe set. In real world, some universe has been given operations, such as the set of natural numbers and the set of real numbers. In applied mathematics, first of all, the real numbers are such objects. Other examples are real-valued functions, the complex numbers and so on. So, it is very natural to ask what would happen if we substitute an algebraic structure for universe set.

For the progression of rough set theory, two most common methods exist: the constructive method and the axiomatic method. In contrast to the constructive approach, the axiomatic approach considers abstract upper and lower approximators subject to certain axioms the primitive notions, and seeks for conditions restraining the axioms to guarantee the existence of a binary relation R on U such that the abstract upper and the abstract lower approximators can be derived from R in the usual way. The axiomatic approach aims to investigate the logical structure of rough sets rather than to develop some methods for applications. The axiomatic approaches can help us to gain much more insight into the logical structure of rough sets and fuzzy rough sets. It also helps unify Pawlak's rough sets and fuzzy rough sets. Using Constructive method, the lower and upper approximations are established from basic concepts like equivalence relations on the universe and neighboring systems. Both lower and upper approximations have equivalence classes as their groundwork. Rough set theory has been integrated with other mathematical theories by implementing these two techniques. In Pawlak rough sets, the equivalence classes are the building blocks for the construction of the lower and upper approximations. The lower approximation of a given set is the union of all the equivalence classes which are the subsets of the set, and the upper approximation is the union of all the equivalence classes which have a non-empty intersection with the set. It is well known that a partition induces an equivalence relation on a set and vice versa. The properties of rough sets can be examined via either partition or equivalence classes.

In rough sets, equivalence classes play an important role in construction of both lower and upper approximations. But sometimes in algebraic structure, as is the case in Γ-semigroups, finding equivalence relation is too difficult. Many authors have worked on this to initiate rough sets without equivalence relations. Some authors substituted an algebraic structure for the universal set and studied the roughness in algebraic structure. Numerous authors have stated about the properties and algebraic structure of rough sets. Combining the theory of rough set with abstract algebra is one of the trends in the theory of rough set. Large number of studies have found a strong link between algebraic systems and rough sets among these research fields. B.Davvaz [8] initiated a generalized rough set or $\mathcal{T}$-rough set with the help of set-valued mapping. The association between rough set and $\mathcal{T}$-rough set is described from this perspective. Specifically, Ideal theory plays a basic role in the study of Γ-semigroups.

Based on constructive method, extensive research has also been carried out to compare the theory of rough sets with other theories such as fuzzy sets, soft sets and conditional events. The relations between rough sets, fuzzy sets, soft sets and algebraic systems have been already considered by many mathematicians. Extensive applications of the rough-fuzzy set theory and rough-soft set theory have been found in various fields. There were many researches investigated the connections and the differences of rough set theory, fuzzy set theory and soft set theory.

## 2     Review of literature

Rough set theory is now rigorous area of research with manifold applications ranging from engineering and computer science to medical diagnosis and social behavior studies. The algebraic approach of rough set was studied by many authors.

Biswas and Nanda [5] had reported earlier rough subgroups. Kuroki [15] made the new proposal for rough ideal concepts in semigroups. With consideration to an ideal of ring, S.Rasouli and B. Davvaz [7,17] raised an idea of rough subring and explored roughness in MV-algebra. Davvaz and Mahdavipour [6] introduced rough modules.

Multiple specialists [14,16,18,19,20,26,27,29,30,32] have focused their research on the congruence relation for rough sets on algebraic structures such as semigroups, groups, rings, commutative rings, ideals, lattices, ternary semigroups and modules. Consequently, the generalized rough set model on algebraic sets can sometimes be utilized without attention of the congruence relation. Davvaz speculated an idea of a set-valued homomorphism of groups to bridge this gap. Initially designed by Davvaz [8], algebraic $\mathcal{T}$-rough sets. $\mathcal{T}$-roughness in semigroups were analyzed by Qiami Xiao [9]. On commutative rings and submodules, S.B.Hosseini [10,11] studied $\mathcal{T}$-rough semiprime ideals. Some aspects of the $\mathcal{T}$-rough ideals and the $\mathcal{T}$-rough fuzzy (prime, primary) ideals on commutative rings were examined by S.B. Hosseini and N. Jafarzadeh [12,13]. In addition to outlining $\mathcal{T}$-rough sets based on lattices, S.B.Hosseini and E.Hosseinpour [28,32] highlighted some of the criteria of a set-valued homomorphism on modules. Ternary semigroups in relation to $\mathcal{T}$-rough ideals were done by Moin A.Ansari and Naveed Yaqoob and Shahida Bashir et.al[14,31].

The algebraic structure of Γ-semigroup was first introduced by Sen and Saha [21-24] to indicate a generalization of semigroup. Later many authors investigated the structure of Γ-semigroups and added some important results related to Γ-semigroups. M.A.Ansari and M.Rais khan [37] were constructed rough (m, n) quasi-Γ-ideals in Γ-semigroups. Anusorn Simuen [33] et al. provided the characteristics of almost bi-Γ-ideals in Γ-semigroups. R.Chinram and N. Yahoob [34-36] reviewed rough fuzzy prime bi-ideals and quasi ideals in semigroups as well as rough prime (m, n) bi-ideals.

In this paper, we have been using set-valued anti-homomorphism to discuss the associations between the lower and upper approximations of some types of (prime)ideals and quotient ideals in Γ-semigroups. We established the concepts of Γ-semigroups and defined the terms $\mathcal{T}$-rough bi-ideal, $\mathcal{T}$-rough quasi-ideal, $\mathcal{T}$-rough interior ideal, $\mathcal{T}$-rough bi-quasi ideal, $\mathcal{T}$-rough bi-interior ideal, $\mathcal{T}$-rough quasi-interior ideal, $\mathcal{T}$-rough bi-quasi-interior ideal, $\mathcal{T}$-rough prime bi-ideal, $\mathcal{T}$-rough prime quasi-ideal, $\mathcal{T}$-rough prime interior ideal, $\mathcal{T}$-rough prime bi-interior ideal, $\mathcal{T}$-rough prime bi-quasi ideal, $\mathcal{T}$-rough prime quasi-interior ideal, $\mathcal{T}$-rough prime bi-quasi-interior ideal and $\mathcal{T}$-rough quotient ideals of Γ-semigroups. we discuss a general mathematical concept, called a Γ-semigroup, which includes all examples and many others as special cases.

The rest of the paper is organized as follows. In section 3, some structural characteristics were described in addition to the definition of Γ-semigroup and the idea of generalised rough approximation operators. In section 4, we propose the definition of the generalised $\mathcal{T}$-rough set on Γ-semigroups and the definition of (strong) set-valued anti-homomorphism on Γ-semigroups with example. In section 5, we have examined some characteristics of the generalized upper and the lower approximations of some types of (prime)ideals and quotient ideals in Γ-semigroups through set-valued anti-homomorphism with examples. In section 6, we have provided the application of rough set theory. The paper is completed with some concluding remarks.

## 3  Preliminaries

In this section, we study Pawlak roughness and generalized roughness in Γ-semigroups.

## Pawlak Approximations in Γ-semigroups

The concept of a rough set was introduced by Pawlak [1]. According to Pawlak, rough set theory is based on the approximations of a set by a pair of sets called lower approximation and upper approximation of that set. Let $\mathcal{U}$ be a non-empty finite set with an equivalence relation $\mathcal{R}$. We say $(\mathcal{U}, \mathcal{R})$ is the approximation space. If $\mathcal{A} \subseteq \mathcal{U}$ can be written as the union of some classes obtained from $\mathcal{R}$, then $\mathcal{A}$ is called definable; otherwise, it is not definable. Therefore, the approximations of $\mathcal{A}$ are as follows:

$$\underline{\mathcal{R}}(\mathcal{A}) = \{x \in \mathcal{U}/[x]_\mathcal{R} \subseteq \mathcal{A}\} \text{ and } \overline{\mathcal{R}}(\mathcal{A}) = \{x \in \mathcal{U}/[x]_\mathcal{R} \cap \mathcal{A} \neq \emptyset\}.$$

The pair $(\underline{\mathcal{R}}(\mathcal{A}), \overline{\mathcal{R}}(\mathcal{A}))$ is a rough set, where $\underline{\mathcal{R}}(\mathcal{A}) \neq \overline{\mathcal{R}}(\mathcal{A})$.

**Definition 3.1.** Let $\rho$ be an equivalence relation on $\mathcal{M}$. Then $\rho$ is called a congruence relation on $\mathcal{M}$ if $(a, b) \in \rho$ implies that $(a\alpha y, b\alpha y) \in \rho$ and $(y\alpha a, y\alpha b) \in \rho$ for all a, b, y∈ $\mathcal{M}$ and α ∈ Γ.

If $\rho$ is a congruence relation on $\mathcal{M}$, then for every x ∈ $\mathcal{M}$, $[x]_\rho$ denotes the congruence class of x with respect to the relation $\rho$.

A congruence relation $\rho$ on $\mathcal{M}$ is called complete, if $[a]_\rho \Gamma [b]_\rho = [a\Gamma b]_\rho$ ∀ a, b∈ $\mathcal{M}$.

**Definition 3.2.** Let $\rho$ be a congruence relation on $\mathcal{M}$. Then the approximation of $\mathcal{M}$ is defined by $\rho(\mathcal{A}) = (\underline{\rho}(\mathcal{A}), \overline{\rho}(\mathcal{A}))$ for every $\mathcal{A} \in \mathcal{P}^*(\mathcal{M})$ is the power set of $\mathcal{M}$, and

$$\underline{\rho}(\mathcal{A}) = \{x \in \mathcal{U}/[x]_\rho \subseteq \mathcal{A}\} \text{ and } \overline{\rho}(\mathcal{A}) = \{x \in \mathcal{U}/[x]_\rho \cap \mathcal{A} \neq \emptyset\}.$$

## Generalized Roughness or $\mathcal{T}$-Roughness in Γ-semigroups

A generalized rough set is the generalization of Pawlak`s rough set. In this case, we use set-valued mappings instead of congruence classes.

**Definition 3.3. [8]** Let $\mathcal{G}$ and $\mathcal{H}$ be two non-empty sets and $\mathcal{B} \subseteq \mathcal{H}$. A set-valued function $\mathcal{T}: \mathcal{G} \to \mathcal{P}^*(\mathcal{H})$ where $\mathcal{P}^*(\mathcal{H})$ represents the set of all non-empty subsets of $\mathcal{H}$. The lower and upper approximations of $\mathcal{B}$ are given by

$$\underline{\mathcal{T}}(\mathcal{B}) = \{x \in \mathcal{G}/\mathcal{T}(x) \subseteq \mathcal{B}\} \text{ and } \overline{\mathcal{T}}(\mathcal{B}) = \{x \in \mathcal{G}/\mathcal{T}(x) \cap \mathcal{B} \neq \emptyset\}$$

The pair $(\underline{\mathcal{T}}(\mathcal{B}), \overline{\mathcal{T}}(\mathcal{B}))$ is known as generalized rough set of $\mathcal{B}$ induced by $\mathcal{T}$ and the generalized operators are denoted by $\underline{\mathcal{T}}$, $\overline{\mathcal{T}}$ respectively.

**Example 1.** Let X = {1,2,3,4}, Y = {a, b, c}. Consider the set-valued mapping $\mathcal{T}: X \to \mathcal{P}^*(Y)$ defined by $\mathcal{T}(1) = \{b\}, \mathcal{T}(2) = \{a, c\}, \mathcal{T}(3) = \{b\}, \mathcal{T}(4) = \{a, b, c\}$. Then,

$\overline{\mathcal{T}}(\{a\}) = \{2,4\}$ $\quad\quad\quad$ $\underline{\mathcal{T}}(\{a\}) = \emptyset$

$\overline{\mathcal{T}}(\{b\}) = \{1,3\}$ $\quad\quad\quad$ $\underline{\mathcal{T}}(\{b\}) = \{1,3\}$

$\overline{\mathcal{T}}(\{c\}) = \{2,4\}$ $\quad\quad\quad$ $\underline{\mathcal{T}}(\{c\}) = \emptyset$

$\overline{\mathcal{T}}(\{a, b\}) = \{1,2,3,4\}$ $\quad\quad\quad$ $\underline{\mathcal{T}}(\{a, b\}) = \{1,3\}$

$\overline{\mathcal{T}}(\{a, c\}) = \{2,4\}$ $\quad\quad\quad$ $\underline{\mathcal{T}}(\{a, c\}) = \{2\}$

$$\overline{\mathcal{T}}(\{b,c\}) = \{1,2,3,4\} \qquad \underline{\mathcal{T}}(\{b,c\}) = \{1,3\}$$

$$\overline{\mathcal{T}}(\{a,b,c\}) = \{1,2,3,4\} \qquad \underline{\mathcal{T}}(\{a,b,c\}) = \{1,2,3,4\}$$

**Proposition 3.4. [9]** Consider two non-empty sets $\mathcal{G}$ and $\mathcal{H}$. Define a set-valued mapping $\mathcal{T}: \mathcal{G} \to \mathcal{P}^*(\mathcal{H})$. Given a generalized approximation space $(\mathcal{G}, \mathcal{H}, \mathcal{T})$. Following are the basic properties of rough approximations, if $\mathcal{X}$ and $\mathcal{Y}$ are non-empty subsets of $\mathcal{H}$:

$$\overline{\mathcal{T}}(\mathcal{X}) = (\underline{\mathcal{T}}(\mathcal{X}^c))^c \qquad \underline{\mathcal{T}}(\mathcal{X}) = (\overline{\mathcal{T}}(\mathcal{X}^c))^c$$

$$\overline{\mathcal{T}}(\mathcal{H}) = \mathcal{G} \qquad \underline{\mathcal{T}}(\mathcal{H}) = \mathcal{G}$$

$$\overline{\mathcal{T}}(\mathcal{X} \cap \mathcal{Y}) \subseteq \overline{\mathcal{T}}(\mathcal{X}) \cap \overline{\mathcal{T}}(\mathcal{Y}) \qquad \underline{\mathcal{T}}(\mathcal{X} \cap \mathcal{Y}) = \underline{\mathcal{T}}(\mathcal{X}) \cap \underline{\mathcal{T}}(\mathcal{Y})$$

$$\overline{\mathcal{T}}(\mathcal{X} \cup \mathcal{Y}) = \overline{\mathcal{T}}(\mathcal{X}) \cup \overline{\mathcal{T}}(\mathcal{Y}) \qquad \underline{\mathcal{T}}(\mathcal{X} \cup \mathcal{Y}) \supseteq \underline{\mathcal{T}}(\mathcal{X}) \cup \underline{\mathcal{T}}(\mathcal{Y})$$

$$\mathcal{X} \subseteq \mathcal{Y} \Rightarrow \overline{\mathcal{T}}(\mathcal{X}) \subseteq \overline{\mathcal{T}}(\mathcal{Y}) \qquad \mathcal{X} \subseteq \mathcal{Y} \Rightarrow \underline{\mathcal{T}}(\mathcal{X}) \subseteq \underline{\mathcal{T}}(\mathcal{Y})$$

**Definition 3.5. [22]** Consider two non-empty sets $\mathcal{M}$ and $\Gamma$. If a mapping $\mathcal{M} \times \Gamma \times \mathcal{M} \to \mathcal{M}$ expressed as $(a, \alpha, b) \to a\alpha b$ satisfied the axiom $(a\alpha b)\beta c = a\alpha(b\beta c)$, $\forall a, b, c \in \mathcal{M}$ and $\alpha, \beta \in \Gamma$, then $\mathcal{M}$ is said to be $\Gamma$-semigroup.

Let $\mathcal{M}$ be a $\Gamma$-semigroup. If $\mathcal{M}_1, \mathcal{M}_2 \subseteq \mathcal{M}$, then we denote $\{a\gamma b \;/\; a \in \mathcal{M}_1, b \in \mathcal{M}_2, \gamma \in \Gamma\}$ by $\mathcal{M}_1 \Gamma \mathcal{M}_2$.

Let $\mathcal{M}$ be a $\Gamma$-semigroup and $\Gamma$ any non-empty set. Let $a\gamma b = ab \; \forall \; a, b \in \mathcal{M}$ and $\gamma \in \Gamma$. It is clear that $\mathcal{M}$ is a $\Gamma$-semigroup. Thus, a semigroup can be considered to be a $\Gamma$-semigroup.

Let $\mathcal{M}$ be a $\Gamma$-semigroup. We define $a.b = a\gamma b \; \forall \; a, b \in \mathcal{M}$ and $\gamma \in \Gamma$. We can show that $(\mathcal{M},.)$ is a semigroup.

**Example 2.** Let $\mathcal{M} = \{1,2,3,4\}$ and $\Gamma = \{\alpha, \beta\}$

| $\alpha$ | 1 | 2 | 3 | 4 |
|---|---|---|---|---|
| 1 | 1 | 3 | 3 | 1 |
| 2 | 3 | 1 | 1 | 3 |
| 3 | 3 | 1 | 1 | 3 |
| 4 | 1 | 3 | 3 | 3 |

| $\beta$ | 1 | 2 | 3 | 4 |
|---|---|---|---|---|
| 1 | 3 | 1 | 1 | 3 |
| 2 | 1 | 3 | 3 | 1 |
| 3 | 1 | 3 | 3 | 1 |
| 4 | 3 | 1 | 1 | 3 |

Here, $(a\alpha b)\beta c = a\alpha(b\beta c)$, $\forall a, b, c \in \mathcal{M}$ and $\alpha, \beta \in \Gamma$. Thus $\mathcal{M}$ is a $\Gamma$-semigroup.

**Definition 3.5.1. [25]** A non-empty subset $\mathcal{S}$ of $\Gamma$-semigroup $\mathcal{M}$ is sub-$\Gamma$-semigroup of $\mathcal{M}$, if $\mathcal{S}\Gamma\mathcal{S} \subseteq \mathcal{S}$.

**Definition 3.5.2. [25]** A non-empty subset $\mathcal{J}$ of $\Gamma$-semigroup $\mathcal{M}$ is called left ideal of $\mathcal{M}$, if $\mathcal{M}\Gamma\mathcal{J} \subseteq \mathcal{J}$ and the subset $\mathcal{J}$ is right ideal of $\mathcal{M}$, if $\mathcal{J}\Gamma\mathcal{M} \subseteq \mathcal{J}$.

**Definition 3.5.3. [25]** Let $\mathcal{M}$ be $\Gamma$-semigroup and $\mathcal{J} \subseteq \mathcal{M}$. Then $\mathcal{J}$ is said to be an ideal of $\mathcal{M}$, if $\mathcal{J}$ is both left ideal and right ideal of $\mathcal{M}$.

**Definition 3.5.4. [38]** A non-empty subset $\mathcal{B}$ of $\Gamma$-semigroup $\mathcal{M}$ is called bi-ideal of $\mathcal{M}$, if $\mathcal{B}$ is sub-$\Gamma$-semigroup of $\mathcal{M}$ and $\mathcal{B}\Gamma\mathcal{M}\Gamma\mathcal{B} \subseteq \mathcal{B}$.

**Definition 3.5.5. [38]** A non-empty subset $Q$ of $\Gamma$-semigroup $\mathcal{M}$ is said to be quasi-ideal of $\mathcal{M}$, if $Q$ is sub-$\Gamma$-semigroup of $\mathcal{M}$ and $Q\Gamma\mathcal{M} \cap \mathcal{M}\Gamma Q \subseteq Q$.

**Definition 3.5.6. [38]** A non-empty subset $\mathcal{J}$ of $\Gamma$-semigroup $\mathcal{M}$ is known as an interior ideal of $\mathcal{M}$, if $\mathcal{J}$ is sub-$\Gamma$-semigroup of $\mathcal{M}$ and $\mathcal{M}\Gamma\mathcal{J}\Gamma\mathcal{M} \subseteq \mathcal{J}$.

**Definition 3.5.7. [38]** A non-empty subset $\mathcal{B}$ of $\Gamma$-semigroup $\mathcal{M}$ is called left bi-quasi ideal of $\mathcal{M}$, if $\mathcal{B}$ is sub-$\Gamma$-semigroup of $\mathcal{M}$ and $\mathcal{M}\Gamma\mathcal{B} \cap \mathcal{B}\Gamma\mathcal{M}\Gamma\mathcal{B} \subseteq \mathcal{B}$ and the subset $\mathcal{B}$ is right bi-quasi ideal of $\mathcal{M}$, if $\mathcal{B}\Gamma\mathcal{M} \cap \mathcal{B}\Gamma\mathcal{M}\Gamma\mathcal{B} \subseteq \mathcal{B}$. Then $\mathcal{B}$ is said to be bi-quasi ideal of $\mathcal{M}$, if $\mathcal{B}$ is both left bi-quasi ideal and right bi-quasi ideal of $\mathcal{M}$.

**Definition 3.5.8. [39]** A non-empty subset $\mathcal{B}$ of $\Gamma$-semigroup $\mathcal{M}$ is called bi-interior ideal of $\mathcal{M}$, if $\mathcal{B}$ is sub-$\Gamma$-semigroup of $\mathcal{M}$ and $\mathcal{M}\Gamma\mathcal{B}\Gamma\mathcal{M} \cap \mathcal{B}\Gamma\mathcal{M}\Gamma\mathcal{B} \subseteq \mathcal{B}$.

**Definition 3.5.9. [39]** A non-empty subset $Q$ of $\Gamma$-semigroup $\mathcal{M}$ is said to be left quasi-interior ideal of $\mathcal{M}$, if $Q$ is sub-$\Gamma$-semigroup of $\mathcal{M}$ and $\mathcal{M}\Gamma Q\Gamma\mathcal{M}\Gamma Q \subseteq Q$ and the subset $Q$ is right quasi-interior ideal of $\mathcal{M}$, if $Q\Gamma\mathcal{M}\Gamma Q\Gamma\mathcal{M} \subseteq Q$. Then $Q$ is said to be quasi-interior ideal of $\mathcal{M}$, if $Q$ is both left quasi-interior ideal and right quasi-interior ideal of $\mathcal{M}$.

**Definition 3.5.10. [39]** A non-empty subset $\mathcal{B}$ of $\Gamma$-semigroup $\mathcal{M}$ is said to be bi-quasi-interior ideal of $\mathcal{M}$, if $\mathcal{B}$ is sub-$\Gamma$-semigroup of $\mathcal{M}$ and $\mathcal{B}\Gamma\mathcal{M}\Gamma\mathcal{B}\Gamma\mathcal{M}\Gamma\mathcal{B} \subseteq \mathcal{B}$.

**Definition 3.5.11. [25]** An ideal $\mathcal{J}$ of $\Gamma$-semigroup $\mathcal{M}$ is prime ideal of $\mathcal{M}$, if $x, y \in \mathcal{M}$, $x\alpha y \in \mathcal{J}$ implies $x \in \mathcal{J}$ or $y \in \mathcal{J}$, $\forall \alpha \in \Gamma$.

**Definition 3.5.12. [35]** A bi-ideal $\mathcal{B}$ of $\Gamma$-semigroup $\mathcal{M}$ is said to be prime bi-ideal of $\mathcal{M}$, if $a, b \in \mathcal{M}$, $a\alpha b \in \mathcal{B}$ implies $a \in \mathcal{B}$ or $b \in \mathcal{B}$, $\forall \alpha \in \Gamma$.

**Definition 3.5.13.** A quasi-ideal $Q$ of $\Gamma$-semigroup $\mathcal{M}$ is prime quasi-ideal of $\mathcal{M}$, if $a, b \in \mathcal{M}$, $a\alpha b \in Q$ implies $a \in Q$ or $b \in Q$, $\forall \alpha \in \Gamma$.

**Definition 3.5.14.** An interior ideal $\mathcal{J}$ of $\Gamma$-semigroup $\mathcal{M}$ is prime interior ideal of $\mathcal{M}$, if $x, y \in \mathcal{M}$, $x\alpha y \in \mathcal{J}$ implies $x \in \mathcal{J}$ or $y \in \mathcal{J}$, $\forall \alpha \in \Gamma$.

**Definition 3.5.15.** A bi-interior ideal $\mathcal{B}$ of $\Gamma$-semigroup $\mathcal{M}$ is said to be prime bi-interior ideal of $\mathcal{M}$, if $a, b \in \mathcal{M}$, $a\alpha b \in \mathcal{B}$ implies $a \in \mathcal{B}$ or $b \in \mathcal{B}$, $\forall \alpha \in \Gamma$.

**Definition 3.5.16.** A (left, right) bi-quasi ideal $\mathcal{B}$ of $\Gamma$-semigroup $\mathcal{M}$ is prime (left, right) bi-quasi ideal of $\mathcal{M}$, if $a, b \in \mathcal{M}$, $a\alpha b \in \mathcal{B}$ implies $a \in \mathcal{B}$ or $b \in \mathcal{B}$, $\forall \alpha \in \Gamma$.

**Definition 3.5.17.** A (left, right) quasi-interior ideal $Q$ of $\Gamma$-semigroup $\mathcal{M}$ is prime (left, right) quasi-interior ideal of $\mathcal{M}$, if $x, y \in \mathcal{M}$, $x\alpha y \in Q$ implies $x \in Q$ or $y \in Q$, $\forall \alpha \in \Gamma$.

**Definition 3.5.18.** A bi-quasi-interior ideal $\mathcal{B}$ of $\Gamma$-semigroup $\mathcal{M}$ is said to be prime bi-quasi-interior ideal of $\mathcal{M}$, if $a, b \in \mathcal{M}$, $a\alpha b \in \mathcal{B}$ implies $a \in \mathcal{B}$ or $b \in \mathcal{B}$, $\forall \alpha \in \Gamma$.

# 4                  Proposed Work

Rough sets were originally proposed in the presence of an equivalence relation. An equivalence relation is sometimes difficult to be obtained in real-world problems due to the vagueness and incompleteness of human knowledge. On the other hand, in applied mathematics we encounter many examples of mathematical objects that can be added or multiplied to each other. From this point of view, we introduce the concept of set-value anti-homomorphism of Γ-semigroups. It's another approach to investigate the properties of generalized approximations by set-valued anti-homomorphism of Γ-semigroups. Moreover, a new algebraic structure called generalized lower and upper approximations of a set of Γ-semigroups with respect to set-valued anti-homomorphism and its properties proved. This paper is an extended idea presented by Davvaz [8].

However, many researchers proposed the rough set which is induced by set-valued homomorphism in groups, rings (with respect to an ideal of rings), commutative rings, modules and investigated some properties of generalized approximations.

**Definition 4.1.** Let $\mathcal{M}_1$ and $\mathcal{M}_2$ be two Γ-semigroups and $\mathcal{G} \subseteq \mathcal{M}_2$ and a set-valued mapping $\mathcal{T}: \mathcal{M}_1 \to \mathcal{P}^*(\mathcal{M}_2)$ where $\mathcal{P}^*(\mathcal{M}_2)$ stands for the collection of all non-empty subsets of $\mathcal{M}_2$. Following are definitions of the generalized lower and upper approximations of $\mathcal{G}$:

$$\underline{\mathcal{T}}(\mathcal{G}) = \{x \in \mathcal{M}_1 / \mathcal{T}(x) \subseteq \mathcal{G}\} \text{ and } \overline{\mathcal{T}}(\mathcal{G}) = \{x \in \mathcal{M}_1 / \mathcal{T}(x) \cap \mathcal{G} \neq \emptyset\}$$

**Definition 4.2.** Mapping $\mathcal{T}$ of $\mathcal{M}_1$ into $\mathcal{P}^*(\mathcal{M}_2)$ that preserves Γ-semigroup operation, $\mathcal{T}(a\alpha b) \supseteq \mathcal{T}(b)\alpha \mathcal{T}(a), \forall a, b \in \mathcal{M}_1$ and $\alpha \in \Gamma$, is known as set-valued anti-homomorphism $\mathcal{T}$ from Γ-semigroup $\mathcal{M}_1$ to Γ-semigroup $\mathcal{M}_2$.

**Example 3.** Let $\mathcal{M} = \{a, b, c\}$ and $\Gamma = \{\alpha\}$ with the following table:

| α | a | b | c |
|---|---|---|---|
| a | a | b | c |
| b | b | b | b |
| c | c | b | b |

Then $\mathcal{M}$ is a Γ-semigroup. Define a set-valued mapping $\mathcal{T}: \mathcal{M} \to \mathcal{P}^*(\mathcal{M})$ by $\mathcal{T}(a) = \{b, c\}, \mathcal{T}(b) = \{a, b, c\}, \mathcal{T}(c) = \{b\}$. Then clearly $\mathcal{T}$ is a set-valued anti-homomorphism.

**Example 4.** Let $\mathcal{M} = \{x_1, x_2, x_3, x_4\}$ and $\Gamma = \{\alpha\}$ with the following table:

| α | $x_1$ | $x_2$ | $x_3$ | $x_4$ |
|---|---|---|---|---|
| $x_1$ | $x_1$ | $x_3$ | $x_3$ | $x_1$ |
| $x_2$ | $x_3$ | $x_1$ | $x_1$ | $x_3$ |
| $x_3$ | $x_3$ | $x_1$ | $x_1$ | $x_3$ |
| $x_4$ | $x_1$ | $x_3$ | $x_3$ | $x_1$ |

Then $\mathcal{M}$ is a Γ-semigroup. Define a set-valued mapping $\mathcal{T}: \mathcal{M} \to \mathcal{P}^*(\mathcal{M})$ by $\mathcal{T}(x_1) = \{x_1, x_2, x_3, x_4\}, \mathcal{T}(x_2) = \{x_1, x_3\}, \mathcal{T}(x_3) = \{x_3\}, \mathcal{T}(x_4) = \{x_4\}$. Then clearly $\mathcal{T}$ is a set-valued anti-homomorphism.

**Definition 4.3.** Consider two Γ-semigroups $\mathcal{M}_1$ and $\mathcal{M}_2$. Define a set-valued mapping $\mathcal{T}: \mathcal{M}_1 \to \mathcal{P}^*(\mathcal{M}_2)$. If $\mathcal{T}(a\alpha b) = \mathcal{T}(b)\alpha \mathcal{T}(a), \forall a, b \in \mathcal{M}_1$ and $\alpha \in \Gamma$, then $\mathcal{T}$ is stated as strong set-valued anti-homomorphism.

**Example 5.** Let $\mathcal{M}_1 = \{1,2\}$ and $\Gamma = \{\alpha\}$ with the following table:

| α | 1 | 2 |
|---|---|---|
| 1 | 1 | 1 |
| 2 | 1 | 2 |

Then $\mathcal{M}_1$ is a $\Gamma$-semigroup and $\mathcal{M}_2 = \{a, b, c\}$ with the following table:

| α | a | b | c |
|---|---|---|---|
| a | a | b | c |
| b | b | b | b |
| c | c | b | b |

Then $\mathcal{M}_2$ is a $\Gamma$-semigroup. Define a set-valued mapping $\mathcal{T}: \mathcal{M}_1 \to \mathcal{P}^*(\mathcal{M}_2)$ Assume $\mathcal{T}(1) = \{c\}, \mathcal{T}(2) = \{a\}$. Clearly $\mathcal{T}$ is a strong set-valued anti-homomorphism.

**Example 6.** Let $\mathcal{M}_1 = \{x, y, z\}$ and $\Gamma = \{\alpha\}$ with the following table:

| α | x | y | z |
|---|---|---|---|
| x | x | x | x |
| y | z | z | z |
| z | x | x | z |

Then $\mathcal{M}_1$ is a $\Gamma$-semigroup and $\mathcal{M}_2 = \{a, b, c\}$ with the following table:

| α | a | b | c |
|---|---|---|---|
| a | a | a | c |
| b | c | c | c |
| c | a | a | a |

Then $\mathcal{M}_2$ is a $\Gamma$-semigroup. Define a set-valued mapping $\mathcal{T}: \mathcal{M}_1 \to \mathcal{P}^*(\mathcal{M}_2)$ Assume $\mathcal{T}(x) = \mathcal{T}(z) = \{c\}, \mathcal{T}(y) = \{b, c\}$. Clearly $\mathcal{T}$ is a strong set-valued anti-homomorphism.

## 5     Results and Discussions

### $\mathcal{T}$-Rough (prime) ideals in $\Gamma$-semigroups

**Theorem 5.1.** Let $\mathcal{M}_1$ and $\mathcal{M}_2$ be two $\Gamma$-semigroups and $\mathcal{T}: \mathcal{M}_1 \to \mathcal{P}^*(\mathcal{M}_2)$ be a set-valued anti-homomorphism. Suppose $\mathcal{A}_1, \mathcal{A}_2 \subseteq \mathcal{M}_2$. Then

(i)   $\overline{\mathcal{T}}(\mathcal{A}_1 \Gamma \mathcal{A}_2) \supseteq \overline{\mathcal{T}}(\mathcal{A}_2) \Gamma \overline{\mathcal{T}}(\mathcal{A}_1)$
(ii)  $\underline{\mathcal{T}}(\mathcal{A}_1 \Gamma \mathcal{A}_2) \supseteq \underline{\mathcal{T}}(\mathcal{A}_2) \underline{\mathcal{T}}(\mathcal{A}_1)$

**Proof.** (i) Assume $p \in \overline{\mathcal{T}}(\mathcal{A}_2) \Gamma \overline{\mathcal{T}}(\mathcal{A}_1)$. Then $p = q\alpha r, \forall \alpha \in \Gamma$ where $q \in \overline{\mathcal{T}}(\mathcal{A}_2), r \in \overline{\mathcal{T}}(\mathcal{A}_1)$. Therefore $\mathcal{T}(q) \cap \mathcal{A}_2 \neq \emptyset$ and $\mathcal{T}(r) \cap \mathcal{A}_1 \neq \emptyset$. There exists $u \in \mathcal{T}(q) \cap \mathcal{A}_2, v \in \mathcal{T}(r) \cap \mathcal{A}_1$. Thus $u\alpha v \in \mathcal{A}_2 \Gamma \mathcal{A}_1$. Since $\mathcal{T}$ is set-valued anti-homomorphism, $u\alpha v \in \mathcal{T}(q)\alpha\mathcal{T}(r) \subseteq \mathcal{T}(r\alpha q) = \mathcal{T}(q\alpha r) = \mathcal{T}(p)$. Hence $u\alpha v \in \mathcal{T}(p) \cap \mathcal{A}_1 \Gamma \mathcal{A}_2$ which imply $p \in \overline{\mathcal{T}}(\mathcal{A}_1 \Gamma \mathcal{A}_2)$. Thus $\overline{\mathcal{T}}(\mathcal{A}_1 \Gamma \mathcal{A}_2) \supseteq \overline{\mathcal{T}}(\mathcal{A}_2) \Gamma \overline{\mathcal{T}}(\mathcal{A}_1)$.

(ii) Let $p \in \underline{\mathcal{T}}(\mathcal{A}_2)\underline{\mathcal{T}}(\mathcal{A}_1)$. Then there exists $q \in \underline{\mathcal{T}}(\mathcal{A}_1), r \in \underline{\mathcal{T}}(\mathcal{A}_2)$ such that $p = r\alpha q, \forall \alpha \in \Gamma$. We have $\mathcal{T}(p) = \mathcal{T}(r\alpha q) = \mathcal{T}(q)\alpha\mathcal{T}(r) \subseteq \mathcal{A}_1 \Gamma \mathcal{A}_2$. So, $p \in \underline{\mathcal{T}}(\mathcal{A}_1 \Gamma \mathcal{A}_2)$. Thus $\underline{\mathcal{T}}(\mathcal{A}_1 \Gamma \mathcal{A}_2) \supseteq \underline{\mathcal{T}}(\mathcal{A}_2) \Gamma \underline{\mathcal{T}}(\mathcal{A}_1)$.

**Example 7.** Consider the $\Gamma$-semigroup $\mathcal{M}$ of Example 3. Define a set-valued mapping $\mathcal{T}: \mathcal{M} \to \mathcal{P}^*(\mathcal{M})$ by $\mathcal{T}(a) = \{b, c\}, \mathcal{T}(b) = \{a, b, c\}, \mathcal{T}(c) = \{b\}$. Then $\mathcal{T}$ is a set-valued anti-homomorphism. Let $\mathcal{A}_1 = \{a\}$ and $\mathcal{A}_2 = \{b\}$. Then $\overline{\mathcal{T}}(\mathcal{A}_1) = \{a, b, c\}$ and $\overline{\mathcal{T}}(\mathcal{A}_2) = \{a, b, c\}$. Therefore $\overline{\mathcal{T}}(\mathcal{A}_2)\Gamma\overline{\mathcal{T}}(\mathcal{A}_1) = \{a, b, c\}$ and $\mathcal{A}_1\Gamma\mathcal{A}_2 = \{b\}$. Thus $\overline{\mathcal{T}}(\mathcal{A}_1\Gamma\mathcal{A}_2) = \{a, b, c\}$. Hence $\overline{\mathcal{T}}(\mathcal{A}_1\Gamma\mathcal{A}_2) \supseteq \overline{\mathcal{T}}(\mathcal{A}_2)\Gamma\overline{\mathcal{T}}(\mathcal{A}_1)$.

**Example 8.** Consider the $\Gamma$-semigroups $\mathcal{M}_1$ and $\mathcal{M}_2$ of Example 6. Define a set-valued mapping $\mathcal{T}:\mathcal{M}_1 \to \mathcal{P}^*(\mathcal{M}_2)$ by $\mathcal{T}(x) = \mathcal{T}(z) = \{c\}$, $\mathcal{T}(y) = \{b, c\}$. Then $\mathcal{T}$ is a strong set-valued anti-homomorphism. Assuming the subsets $\{a\},\{c\},\{b, c\}$ and let $\mathcal{A}_1 = \{a, b, c\}$ and $\mathcal{A}_2 = \{b, c\}$. Then $\underline{\mathcal{T}}(\mathcal{A}_1) = \{a, b, c\}$ and $\underline{\mathcal{T}}(\mathcal{A}_2) = \{b, c\}$. Therefore $\underline{\mathcal{T}}(\mathcal{A}_2)\Gamma\underline{\mathcal{T}}(\mathcal{A}_1) = \{a, c\}$ and $\mathcal{A}_1\Gamma\mathcal{A}_2 = \{a, c\}$. Thus $\underline{\mathcal{T}}(\mathcal{A}_1\Gamma\mathcal{A}_2) = \{a, c\}$. Hence $\underline{\mathcal{T}}(\mathcal{A}_1\Gamma\mathcal{A}_2) \supseteq \underline{\mathcal{T}}(\mathcal{A}_2)\Gamma\underline{\mathcal{T}}(\mathcal{A}_1)$.

**Theorem 5.2.** Consider $\mathcal{M}_1$ and $\mathcal{M}_2$ are two $\Gamma$-semigroups and define $\mathcal{T}:\mathcal{M}_1 \to \mathcal{P}^*(\mathcal{M}_2)$ set-valued anti-homomorphism. If $\mathcal{S}$ is a sub-$\Gamma$-semigroup of $\mathcal{M}_2$ then

(i) $\overline{\mathcal{T}}(\mathcal{S})$ is a sub-$\Gamma$-semigroup of $\mathcal{M}_1$, if $\overline{\mathcal{T}}(\mathcal{S}) \neq \emptyset$.
(ii) $\underline{\mathcal{T}}(\mathcal{S})$ is a sub-$\Gamma$-semigroup of $\mathcal{M}_1$, if $\mathcal{T}$ is strong set-valued anti-homomorphism.

**Proof.** (i) Since $\mathcal{S}$ is a sub-$\Gamma$-semigroup of $\mathcal{M}_2$, $\mathcal{S}\Gamma\mathcal{S} \subseteq \mathcal{S}$. By proposition 3.2. and theorem 5.1., $\overline{\mathcal{T}}(\mathcal{S})\Gamma\overline{\mathcal{T}}(\mathcal{S}) \subseteq \overline{\mathcal{T}}(\mathcal{S}\Gamma\mathcal{S}) \subseteq \overline{\mathcal{T}}(\mathcal{S})$. Hence $\overline{\mathcal{T}}(\mathcal{S})$ is a sub-$\Gamma$-semigroup of $\mathcal{M}_1$.

(ii) Let $\mathcal{S}$ be a sub-$\Gamma$-semigroup of $\mathcal{M}_2$. By definition, $\mathcal{S}\Gamma\mathcal{S} \subseteq \mathcal{S}$. By proposition 3.2 and theorem 5.1, $\underline{\mathcal{T}}(\mathcal{S})\Gamma\underline{\mathcal{T}}(\mathcal{S}) \subseteq \underline{\mathcal{T}}(\mathcal{S}\Gamma\mathcal{S}) \subseteq \underline{\mathcal{T}}(\mathcal{S})$. Thus $\underline{\mathcal{T}}(\mathcal{S})$ is a sub-$\Gamma$-semigroup of $\mathcal{M}_1$.

**Example 9.** Consider the $\Gamma$-semigroup $\mathcal{M}$ of Example 3. Define a set-valued mapping $\mathcal{T}:\mathcal{M} \to \mathcal{P}^*(\mathcal{M})$ by $\mathcal{T}(x_1) = \{x_1, x_2, x_3, x_4\}, \mathcal{T}(x_2) = \{x_1, x_3\}, \mathcal{T}(x_3) = \{x_3\}, \mathcal{T}(x_4) = \{x_4\}$. Then clearly $\mathcal{T}$ is a set-valued anti-homomorphism. Let $\mathcal{S} = \{x_1, x_2, x_3\}$ be a subset of $\mathcal{M}$. Assuming the subsets $\{x_1, x_3\}, \{x_2, x_4\}$. Then $\mathcal{S}$ is a sub-$\Gamma$-semigroup. $\overline{\mathcal{T}}(\mathcal{S}) = \{a, b, c, d\}, \mathcal{S}\Gamma\mathcal{S} = \{a, c\}$ and $\overline{\mathcal{T}}(\mathcal{S}\Gamma\mathcal{S}) = \{a, c\}$. Therefore, $\overline{\mathcal{T}}(\mathcal{S}\Gamma\mathcal{S}) \subseteq \overline{\mathcal{T}}(\mathcal{S})$. Thus $\overline{\mathcal{T}}(\mathcal{S})$ is a sub-$\Gamma$-semigroup of $\mathcal{M}$.

**Example 10.** Consider the $\Gamma$-semigroups $\mathcal{M}_1$ and $\mathcal{M}_2$ of Example 5. Define a set-valued mapping $\mathcal{T}:\mathcal{M}_1 \to \mathcal{P}^*(\mathcal{M}_2)$ Assume $\mathcal{T}(1) = \{c\}, \mathcal{T}(2) = \{a\}$. Then $\mathcal{T}$ is a strong set-valued anti-homomorphism. Let $\mathcal{S} = \{b, c\} \subset \mathcal{M}_2$. Assuming the subsets $\{b\}, \{c\}, \{a, c\}$. Then $\mathcal{S}$ is a sub-$\Gamma$-semigroup. $\underline{\mathcal{T}}(\mathcal{S}) = \{b, c\}, \mathcal{S}\Gamma\mathcal{S} = \{b\}$ and $\underline{\mathcal{T}}(\mathcal{S}\Gamma\mathcal{S}) = \{b\}$. Therefore $\underline{\mathcal{T}}(\mathcal{S}\Gamma\mathcal{S}) \subseteq \underline{\mathcal{T}}(\mathcal{S})$. Thus $\underline{\mathcal{T}}(\mathcal{S})$ is a sub-$\Gamma$-semigroup of $\mathcal{M}_1$.

**Theorem 5.3.** Consider $\Gamma$-semigroups $\mathcal{M}_1$ and $\mathcal{M}_2$. Defining $\mathcal{T}:\mathcal{M}_1 \to \mathcal{P}^*(\mathcal{M}_2)$ a set-valued anti-homomorphism. Let $\mathcal{I}$ be a left ideal of $\mathcal{M}_2$. Then

(i) $\overline{\mathcal{T}}(\mathcal{I})$ is a left ideal of $\mathcal{M}_1$, if $\overline{\mathcal{T}}(\mathcal{I}) \neq \emptyset$.
(ii) $\underline{\mathcal{T}}(\mathcal{I})$ is a left ideal of $\mathcal{M}_1$, if $\underline{\mathcal{T}}(\mathcal{I}) \neq \emptyset$ and $\mathcal{T}$ is strong set-valued anti-homomorphism.

**Proof.** (i) Given $\mathcal{I}$ is a left ideal of $\mathcal{M}_2$. Then $\mathcal{M}_2\Gamma\mathcal{I} \subseteq \mathcal{I}$. By proposition, $\overline{\mathcal{T}}(\mathcal{M}_2) = \mathcal{M}_1$. We have, $\mathcal{M}_1\Gamma\overline{\mathcal{T}}(\mathcal{I}) = \overline{\mathcal{T}}(\mathcal{M}_2)\Gamma\overline{\mathcal{T}}(\mathcal{I}) = \overline{\mathcal{T}}(\mathcal{I})\Gamma\overline{\mathcal{T}}(\mathcal{M}_2) \subseteq \overline{\mathcal{T}}(\mathcal{M}_2\Gamma\mathcal{I}) \subseteq \overline{\mathcal{T}}(\mathcal{I})$. Thus $\overline{\mathcal{T}}(\mathcal{I})$ is a left ideal of $\mathcal{M}_1$.

(ii) Since $\mathcal{I}$ is a left ideal of $\mathcal{M}_2$, we have $\mathcal{M}_2\Gamma\mathcal{I} \subseteq \mathcal{I}$. So, $\underline{\mathcal{T}}(\mathcal{M}_2\Gamma\mathcal{I}) \subseteq \underline{\mathcal{T}}(\mathcal{I})$. we have $\underline{\mathcal{T}}(\mathcal{M}_2) = \mathcal{M}_1$. Then $\mathcal{M}_1\Gamma\underline{\mathcal{T}}(\mathcal{I}) = \underline{\mathcal{T}}(\mathcal{M}_2)\Gamma\underline{\mathcal{T}}(\mathcal{I}) = \underline{\mathcal{T}}(\mathcal{I})\Gamma\underline{\mathcal{T}}(\mathcal{M}_2) \subseteq \underline{\mathcal{T}}(\mathcal{M}_2\Gamma\mathcal{I}) \subseteq \underline{\mathcal{T}}(\mathcal{I})$. Thus $\underline{\mathcal{T}}(\mathcal{I})$ is a left ideal of $\mathcal{M}_1$, if it is non-empty.

**Example 11.** Consider the $\Gamma$-semigroups $\mathcal{M}_1$ and $\mathcal{M}_2$ of Example 4. Define a set-valued mapping $\mathcal{T}: \mathcal{M} \to \mathcal{P}^*(\mathcal{M})$ by $\mathcal{T}(x_1) = \{x_1, x_2, x_3, x_4\}, \mathcal{T}(x_2) = \{x_1, x_3\}, \mathcal{T}(x_3) = \{x_3\}, \mathcal{T}(x_4) = \{x_4\}$. Then $\mathcal{T}$ is a set-valued anti-homomorphism. Let $\mathcal{I} = \{x_1, x_3, x_4\} \subset \mathcal{M}_2$. Then $\mathcal{I}$ is a left ideal of $\mathcal{M}_2$. $\overline{\mathcal{T}}(\mathcal{I}) = \{x_1, x_2, x_3, x_4\}$, $\mathcal{M}\Gamma\mathcal{I} = \{x_1, x_3\}$, $\overline{\mathcal{T}}(\mathcal{M}\Gamma\mathcal{I}) = \{x_1, x_2, x_3, x_4\}$. Therefore $\overline{\mathcal{T}}(\mathcal{M}\Gamma\mathcal{I}) \subseteq \overline{\mathcal{T}}(\mathcal{I})$. Thus $\overline{\mathcal{T}}(\mathcal{I})$ is a left ideal of $\mathcal{M}$.

**Example 12.** Consider the $\Gamma$-semigroups $\mathcal{M}_1$ and $\mathcal{M}_2$ of Example 5. Define a set-valued mapping $\mathcal{T}: \mathcal{M}_1 \to \mathcal{P}^*(\mathcal{M}_2)$ Assume $\mathcal{T}(1) = \{c\}, \mathcal{T}(2) = \{a\}$. Then $\mathcal{T}$ is a strong set-valued anti-homomorphism. Let $\mathcal{I} = \{b, c\} \subset \mathcal{M}_2$. Then $\mathcal{I}$ is a left ideal of $\mathcal{M}_2$. $\underline{\mathcal{T}}(\mathcal{I}) = \{c\}, \mathcal{M}_2\Gamma\mathcal{I} = \{b, c\}$ and $\underline{\mathcal{T}}(\mathcal{M}_2\Gamma\mathcal{I}) = \{c\}$. Hence $\underline{\mathcal{T}}(\mathcal{M}_2\Gamma\mathcal{I}) \subseteq \underline{\mathcal{T}}(\mathcal{I})$. Thus $\underline{\mathcal{T}}(\mathcal{I})$ is a left ideal of $\mathcal{M}_1$.

**Corollary 1.** Let $\mathcal{T}: \mathcal{M}_1 \to \mathcal{P}^*(\mathcal{M}_2)$ be a set-valued anti-homomorphism. If $\mathcal{I}$ is a right ideal of $\mathcal{M}_2$. Then

(i) $\overline{\mathcal{T}}(\mathcal{I})$ is a right ideal of $\mathcal{M}_1$, if $\overline{\mathcal{T}}(\mathcal{I}) \neq \emptyset$.
(ii) $\underline{\mathcal{T}}(\mathcal{I})$ is a right ideal of $\mathcal{M}_1$, if $\underline{\mathcal{T}}(\mathcal{I}) \neq \emptyset$ and $\mathcal{T}$ is strong set-valued anti-homomorphism.

**Corollary 2.** Let $\mathcal{T}: \mathcal{M}_1 \to \mathcal{P}^*(\mathcal{M}_2)$ be a set-valued anti-homomorphism. If $\mathcal{I}$ is an ideal of $\mathcal{M}_2$. Then

(i) $\overline{\mathcal{T}}(\mathcal{I})$ is an ideal of $\mathcal{M}_1$, if $\overline{\mathcal{T}}(\mathcal{I}) \neq \emptyset$.
(ii) $\underline{\mathcal{T}}(\mathcal{I})$ is an ideal of $\mathcal{M}_1$, if $\underline{\mathcal{T}}(\mathcal{I}) \neq \emptyset$ and $\mathcal{T}$ is strong set-valued anti-homomorphism.

**Theorem 5.4.** Define $\mathcal{T}: \mathcal{M}_1 \to \mathcal{P}^*(\mathcal{M}_2)$ a set-valued anti-homomorphism where $\mathcal{M}_1$ and $\mathcal{M}_2$ are two $\Gamma$-semigroups. Given $\mathcal{B}$ is a bi-ideal of $\mathcal{M}_2$. Then

(i) $\overline{\mathcal{T}}(\mathcal{B})$ is a bi-ideal of $\mathcal{M}_1$, if $\mathcal{T}(x) \neq \emptyset, \forall x \in \mathcal{M}_1$ and $\overline{\mathcal{T}}(\mathcal{B}) \neq \emptyset$.
(ii) $\underline{\mathcal{T}}(\mathcal{B})$ is a bi-ideal of $\mathcal{M}_1$, if $\underline{\mathcal{T}}(\mathcal{B}) \neq \emptyset$ and $\mathcal{T}$ is strong set-valued anti-homomorphism.

**Proof.** (i) Let $\mathcal{B}$ be a bi-ideal of $\mathcal{M}_2$. Then $\mathcal{B}\Gamma\mathcal{M}_2\Gamma\mathcal{B} \subseteq \mathcal{B}$. By proposition $\overline{\mathcal{T}}(\mathcal{M}_2) = \mathcal{M}_1$. We have, $\overline{\mathcal{T}}(\mathcal{B})\Gamma\mathcal{M}_1\Gamma\overline{\mathcal{T}}(\mathcal{B}) = \overline{\mathcal{T}}(\mathcal{B})\Gamma\overline{\mathcal{T}}(\mathcal{M}_2)\Gamma\overline{\mathcal{T}}(\mathcal{B}) \subseteq \overline{\mathcal{T}}(\mathcal{B}\Gamma\mathcal{M}_2\Gamma\mathcal{B}) \subseteq \overline{\mathcal{T}}(\mathcal{B})$. Therefore $\overline{\mathcal{T}}(\mathcal{B})$ is a bi-ideal of $\mathcal{M}_1$.

(ii) Since $\mathcal{B}$ is a bi-ideal of $\mathcal{M}_2$, we have $\mathcal{B}\Gamma\mathcal{M}_2\Gamma\mathcal{B} \subseteq \mathcal{B}$. Then $\underline{\mathcal{T}}(\mathcal{B}\Gamma\mathcal{M}_2\Gamma\mathcal{B}) \subseteq \underline{\mathcal{T}}(\mathcal{B})$. By proposition 3.2, we have $\underline{\mathcal{T}}(\mathcal{M}_2) = \mathcal{M}_1$. Then $\underline{\mathcal{T}}(\mathcal{B})\Gamma\mathcal{M}_1\Gamma\underline{\mathcal{T}}(\mathcal{B}) = \underline{\mathcal{T}}(\mathcal{B})\Gamma\underline{\mathcal{T}}(\mathcal{M}_2)\Gamma\underline{\mathcal{T}}(\mathcal{B}) \subseteq \underline{\mathcal{T}}(\mathcal{B}\Gamma\mathcal{M}_2\Gamma\mathcal{B}) \subseteq \underline{\mathcal{T}}(\mathcal{B})$. Thus $\underline{\mathcal{T}}(\mathcal{B})$ is a bi-ideal of $\mathcal{M}_1$.

**Theorem 5.5.** The mapping $\mathcal{T}: \mathcal{M}_1 \to \mathcal{P}^*(\mathcal{M}_2)$ is a set-valued anti-homomorphism where $\mathcal{M}_1$ and $\mathcal{M}_2$ represent $\Gamma$-semigroups. Let $\mathcal{I}$ be an interior ideal of $\mathcal{M}_2$. Then

(i) $\overline{\mathcal{T}}(\mathcal{I})$ is an interior ideal of $\mathcal{M}_1$, if $\overline{\mathcal{T}}(\mathcal{I}) \neq \emptyset$
(ii) $\underline{\mathcal{T}}(\mathcal{I})$ is an interior ideal of $\mathcal{M}_1$, if $\underline{\mathcal{T}}(\mathcal{I}) \neq \emptyset$ and $\mathcal{T}$ is strong set-valued anti-homomorphism.

**Proof.** (i) Since $\mathcal{J}$ is an interior ideal of $\mathcal{M}_2$, $\mathcal{M}_2\Gamma\mathcal{J}\Gamma\mathcal{M}_2 \subseteq \mathcal{J}$. By proposition $\overline{\mathcal{T}}(\mathcal{M}_2) = \mathcal{M}_1$. We have, $\mathcal{M}_1\Gamma\overline{\mathcal{T}}(\mathcal{J})\Gamma\mathcal{M}_1 = \overline{\mathcal{T}}(\mathcal{M}_2)\Gamma\overline{\mathcal{T}}(\mathcal{J})\Gamma\overline{\mathcal{T}}(\mathcal{M}_2) \subseteq \overline{\mathcal{T}}(\mathcal{M}_2\Gamma\mathcal{J}\Gamma\mathcal{M}_2) \subseteq \overline{\mathcal{T}}(\mathcal{J})$. Hence $\overline{\mathcal{T}}(\mathcal{J})$ is an interior ideal of $\mathcal{M}_1$.

(ii) Let $\mathcal{J}$ be an interior ideal of $\mathcal{M}_2$. By definition, $\mathcal{M}_2\Gamma\mathcal{J}\Gamma\mathcal{M}_2 \subseteq \mathcal{J}$. By property, $\underline{\mathcal{T}}(\mathcal{M}_2) = \mathcal{M}_1$. Now, $\underline{\mathcal{T}}(\mathcal{M}_2\Gamma\mathcal{J}\Gamma\mathcal{M}_2) \subseteq \underline{\mathcal{T}}(\mathcal{J}_1)$. Then, $\mathcal{M}_1\Gamma\underline{\mathcal{T}}(\mathcal{J})\Gamma\mathcal{M}_1 = \underline{\mathcal{T}}(\mathcal{M}_2)\Gamma\underline{\mathcal{T}}(\mathcal{J})\Gamma\underline{\mathcal{T}}(\mathcal{M}_2) \subseteq \underline{\mathcal{T}}(\mathcal{M}_2\Gamma\mathcal{J}\Gamma\mathcal{M}_2) \subseteq \underline{\mathcal{T}}(\mathcal{J})$. Thus $\underline{\mathcal{T}}(\mathcal{J})$ is an interior ideal of $\mathcal{M}_1$.

**Theorem 5.6.** Let $\mathcal{M}_1$ and $\mathcal{M}_2$ be two $\Gamma$-semigroups and $\mathcal{T}: \mathcal{M}_1 \to \mathcal{P}^*(\mathcal{M}_2)$ be a set-valued anti-homomorphism. If $\mathcal{Q}$ is a quasi-ideal of $\mathcal{M}_2$, then

(i) $\overline{\mathcal{T}}(\mathcal{Q})$ is a quasi-ideal of $\mathcal{M}_1$, if $\mathcal{T}(x) \neq \emptyset, \forall\, x \in \mathcal{M}_1$ and $\overline{\mathcal{T}}(\mathcal{Q}) \neq \emptyset$.
(ii) $\underline{\mathcal{T}}(\mathcal{Q})$ is a quasi-ideal of $\mathcal{M}_1$, if $\underline{\mathcal{T}}(\mathcal{Q}) \neq \emptyset$ and $\mathcal{T}$ is strong set-valued anti-homomorphism.

**Proof.** (i) Since $\mathcal{Q}$ is a quasi-ideal of $\mathcal{M}_2$, $\mathcal{Q}\Gamma\mathcal{M}_2 \cap \mathcal{M}_2\Gamma\mathcal{Q} \subseteq \mathcal{Q}$. By proposition $\overline{\mathcal{T}}(\mathcal{M}_2) = \mathcal{M}_1$. We have, $\overline{\mathcal{T}}(\mathcal{Q}\Gamma\mathcal{M}_2 \cap \mathcal{M}_2\Gamma\mathcal{Q}) \subseteq \overline{\mathcal{T}}(\mathcal{Q}\Gamma\mathcal{M}_2) \cap \overline{\mathcal{T}}(\mathcal{M}_2\Gamma\mathcal{Q}) \subseteq \overline{\mathcal{T}}(\mathcal{M}_2)\Gamma\overline{\mathcal{T}}(\mathcal{Q}) \cap \overline{\mathcal{T}}(\mathcal{Q})\Gamma\overline{\mathcal{T}}(\mathcal{M}_2) = \mathcal{M}_1\Gamma\overline{\mathcal{T}}(\mathcal{Q}) \cap \overline{\mathcal{T}}(\mathcal{Q})\Gamma\mathcal{M}_1 \subseteq \overline{\mathcal{T}}(\mathcal{Q})$. Therefore $\overline{\mathcal{T}}(\mathcal{Q})$ is a quasi-ideal of $\mathcal{M}_1$.

(ii) Given $\mathcal{Q}$ is a quasi-ideal of $\mathcal{M}_2$, $\mathcal{Q}\Gamma\mathcal{M}_2 \cap \mathcal{M}_2\Gamma\mathcal{Q} \subseteq \mathcal{Q}$. By proposition $\overline{\mathcal{T}}(\mathcal{M}_2) = \mathcal{M}_1$. Now, $(\mathcal{Q}\Gamma\mathcal{M}_2 \cap \mathcal{M}_2\Gamma\mathcal{Q}) = \underline{\mathcal{T}}(\mathcal{Q}\Gamma\mathcal{M}_2) \cap \underline{\mathcal{T}}(\mathcal{M}_2\Gamma\mathcal{Q}) \subseteq \underline{\mathcal{T}}(\mathcal{M}_2)\Gamma\underline{\mathcal{T}}(\mathcal{Q}) \cap \underline{\mathcal{T}}(\mathcal{Q})\Gamma\underline{\mathcal{T}}(\mathcal{M}_2) = \mathcal{M}_1\Gamma\,\underline{\mathcal{T}}(\mathcal{Q}) \cap \underline{\mathcal{T}}(\mathcal{Q})\Gamma\mathcal{M}_1 \subseteq \underline{\mathcal{T}}(\mathcal{Q})$. Thus $\underline{\mathcal{T}}(\mathcal{Q}\Gamma\mathcal{M}_2 \cap \mathcal{M}_2\Gamma\mathcal{Q}) \subseteq \underline{\mathcal{T}}(\mathcal{Q})$. Hence $\underline{\mathcal{T}}(\mathcal{Q})$ is a quasi-ideal of $\mathcal{M}_1$.

**Theorem 5.7.** Let $\mathcal{M}_1$ and $\mathcal{M}_2$ be two $\Gamma$-semigroups and $\mathcal{T}: \mathcal{M}_1 \to \mathcal{P}^*(\mathcal{M}_2)$ be a set-valued anti-homomorphism. If $\mathcal{B}$ is a bi-interior ideal of $\mathcal{M}_2$, then

(i) $\overline{\mathcal{T}}(\mathcal{B})$ is a bi-interior ideal of $\mathcal{M}_1$, if $\mathcal{T}(x) \neq \emptyset, \forall\, x \in \mathcal{M}_1$ and $\overline{\mathcal{T}}(\mathcal{B}) \neq \emptyset$.
(ii) $\underline{\mathcal{T}}(\mathcal{B})$ is a bi-interior ideal of $\mathcal{M}_1$, if $\underline{\mathcal{T}}(\mathcal{B}) \neq \emptyset$ and $\mathcal{T}$ is strong set-valued anti-homomorphism.

**Proof.** (i) Since $\mathcal{B}$ is a bi-interior ideal of $\mathcal{M}_2$, $\mathcal{M}_2\Gamma\mathcal{B}\Gamma\mathcal{M}_2 \cap \mathcal{B}\Gamma\mathcal{M}_2\Gamma\mathcal{B} \subseteq \mathcal{B}$. By proposition $\overline{\mathcal{T}}(\mathcal{M}_2) = \mathcal{M}_1$. We have, $\overline{\mathcal{T}}(\mathcal{M}_2\Gamma\mathcal{B}\Gamma\mathcal{M}_2 \cap \mathcal{B}\Gamma\mathcal{M}_2\Gamma\mathcal{B}) \subseteq \overline{\mathcal{T}}(\mathcal{M}_2\Gamma\mathcal{B}\Gamma\mathcal{M}_2) \cap \overline{\mathcal{T}}(\mathcal{B}\Gamma\mathcal{M}_2\Gamma\mathcal{B}) \subseteq \overline{\mathcal{T}}(\mathcal{M}_2)\Gamma\overline{\mathcal{T}}(\mathcal{B})\Gamma\overline{\mathcal{T}}(\mathcal{M}_2) \cap \overline{\mathcal{T}}(\mathcal{B})\Gamma\overline{\mathcal{T}}(\mathcal{M}_2)\Gamma\overline{\mathcal{T}}(\mathcal{B}) = \mathcal{M}_1\Gamma\overline{\mathcal{T}}(\mathcal{B})\Gamma\mathcal{M}_1 \cap \overline{\mathcal{T}}(\mathcal{B})\Gamma\mathcal{M}_1\Gamma\overline{\mathcal{T}}(\mathcal{B}) = \mathcal{M}_1\Gamma\overline{\mathcal{T}}(\mathcal{B}) \subseteq \overline{\mathcal{T}}(\mathcal{B})$. Therefore $\overline{\mathcal{T}}(\mathcal{B})$ is a bi-interior ideal of $\mathcal{M}_1$.

(ii) Given $\mathcal{B}$ is a bi-interior ideal of $\mathcal{M}_2$, $\mathcal{M}_2\Gamma\mathcal{B}\Gamma\mathcal{M}_2 \cap \mathcal{B}\Gamma\mathcal{M}_2\Gamma\mathcal{B} \subseteq \mathcal{B}$. By proposition $\underline{\mathcal{T}}(\mathcal{M}_2) = \mathcal{M}_1$. Now, $\underline{\mathcal{T}}(\mathcal{M}_2\Gamma\mathcal{B}\Gamma\mathcal{M}_2 \cap \mathcal{B}\Gamma\mathcal{M}_2\Gamma\mathcal{B}) = \underline{\mathcal{T}}(\mathcal{M}_2\Gamma\mathcal{B}\Gamma\mathcal{M}_2) \cap \underline{\mathcal{T}}(\mathcal{B}\Gamma\mathcal{M}_2\Gamma\mathcal{B}) \subseteq \underline{\mathcal{T}}(\mathcal{M}_2)\Gamma\underline{\mathcal{T}}(\mathcal{B})\Gamma\underline{\mathcal{T}}(\mathcal{M}_2) \cap \underline{\mathcal{T}}(\mathcal{B})\Gamma\underline{\mathcal{T}}(\mathcal{M}_2)\Gamma\underline{\mathcal{T}}(\mathcal{B}) = \mathcal{M}_1\Gamma\,\underline{\mathcal{T}}(\mathcal{B})\Gamma\mathcal{M}_1 \cap \underline{\mathcal{T}}(\mathcal{B})\Gamma\mathcal{M}_1\Gamma\underline{\mathcal{T}}(\mathcal{B}) \subseteq \underline{\mathcal{T}}(\mathcal{B})$. Thus $\underline{\mathcal{T}}(\mathcal{M}_2\Gamma\mathcal{B}\Gamma\mathcal{M}_2 \cap \mathcal{B}\Gamma\mathcal{M}_2\Gamma\mathcal{B}) \subseteq \underline{\mathcal{T}}(\mathcal{B})$. Hence $\underline{\mathcal{T}}(\mathcal{B})$ is a bi-interior ideal of $\mathcal{M}_1$.

**Theorem 5.8.** Let $\mathcal{M}_1$ and $\mathcal{M}_2$ be two $\Gamma$-semigroups and $\mathcal{T}: \mathcal{M}_1 \to \mathcal{P}^*(\mathcal{M}_2)$ be a set-valued anti-homomorphism. If $\mathcal{B}$ is a left bi-quasi ideal of $\mathcal{M}_2$, then

(i) $\overline{\mathcal{T}}(\mathcal{B})$ is a left bi-quasi ideal of $\mathcal{M}_1$, if $\mathcal{T}(x) \neq \emptyset, \forall\, x \in \mathcal{M}_1$ and $\overline{\mathcal{T}}(\mathcal{B}) \neq \emptyset$.
(ii) $\underline{\mathcal{T}}(\mathcal{B})$ is a left bi-quasi ideal of $\mathcal{M}_1$, if $\underline{\mathcal{T}}(\mathcal{B}) \neq \emptyset$ and $\mathcal{T}$ is strong set-valued anti-homomorphism.

**Proof.** (i) Since $\mathcal{B}$ is a left bi-quasi ideal of $\mathcal{M}_2$, $\mathcal{M}_2\Gamma\mathcal{B} \cap \mathcal{B}\Gamma\mathcal{M}_2\Gamma\mathcal{B} \subseteq \mathcal{B}$. By proposition $\overline{\mathcal{T}}(\mathcal{M}_2) = \mathcal{M}_1$. We have, $\overline{\mathcal{T}}(\mathcal{M}_2\Gamma\mathcal{B} \cap \mathcal{B}\Gamma\mathcal{M}_2\Gamma\mathcal{B}) \subseteq \overline{\mathcal{T}}(\mathcal{M}_2\Gamma\mathcal{B}) \cap \overline{\mathcal{T}}(\mathcal{B}\Gamma\mathcal{M}_2\Gamma\mathcal{B}) \subseteq \overline{\mathcal{T}}(\mathcal{B})\Gamma\overline{\mathcal{T}}(\mathcal{M}_2) \cap \overline{\mathcal{T}}(\mathcal{B})\Gamma\overline{\mathcal{T}}(\mathcal{M}_2)\Gamma\overline{\mathcal{T}}(\mathcal{B}) = \overline{\mathcal{T}}(\mathcal{B})\Gamma\mathcal{M}_1 \cap \overline{\mathcal{T}}(\mathcal{B})\Gamma\mathcal{M}_1\Gamma\overline{\mathcal{T}}(\mathcal{B}) = \mathcal{M}_1\Gamma\overline{\mathcal{T}}(\mathcal{B}) = \overline{\mathcal{T}}(\mathcal{B})$. Therefore $\overline{\mathcal{T}}(\mathcal{B})$ is a left bi-quasi ideal of $\mathcal{M}_1$.

(ii) Given $\mathcal{B}$ is a left bi-quasi ideal of $\mathcal{M}_2$, $\mathcal{M}_2\Gamma\mathcal{B} \cap \mathcal{B}\Gamma\mathcal{M}_2\Gamma\mathcal{B} \subseteq \mathcal{B}$ By proposition $\underline{\mathcal{T}}(\mathcal{M}_2) = \mathcal{M}_1$. Now, $\underline{\mathcal{T}}(\mathcal{M}_2\Gamma\mathcal{B} \cap \mathcal{B}\Gamma\mathcal{M}_2\Gamma\mathcal{B}) = \underline{\mathcal{T}}(\mathcal{M}_2\Gamma\mathcal{B}) \cap \underline{\mathcal{T}}(\mathcal{B}\Gamma\mathcal{M}_2\Gamma\mathcal{B}) \subseteq \underline{\mathcal{T}}(\mathcal{B})\Gamma\underline{\mathcal{T}}(\mathcal{M}_2) \cap \underline{\mathcal{T}}(\mathcal{B})\Gamma\underline{\mathcal{T}}(\mathcal{M}_2)\Gamma\underline{\mathcal{T}}(\mathcal{B}) = \underline{\mathcal{T}}(\mathcal{B})\Gamma\mathcal{M}_1 \cap \underline{\mathcal{T}}(\mathcal{B})\Gamma\mathcal{M}_1\Gamma\underline{\mathcal{T}}(\mathcal{B}) = \underline{\mathcal{T}}(\mathcal{B})\Gamma\mathcal{M}_1 = \underline{\mathcal{T}}(\mathcal{B})$. Thus $\underline{\mathcal{T}}(\mathcal{M}_2\Gamma\mathcal{B} \cap \mathcal{B}\Gamma\mathcal{M}_2\Gamma\mathcal{B}) \subseteq \underline{\mathcal{T}}(\mathcal{B})$. Hence $\underline{\mathcal{T}}(\mathcal{B})$ is a left bi-quasi ideal of $\mathcal{M}_1$.

**Corollary 3.** Let $\mathcal{T}: \mathcal{M}_1 \to \mathcal{P}^*(\mathcal{M}_2)$ be a set-valued anti-homomorphism. If $\mathcal{B}$ is a right bi-quasi ideal of $\mathcal{M}_2$. Then

(i) $\overline{\mathcal{T}}(\mathcal{B})$ is a right bi-quasi ideal of $\mathcal{M}_1$, if $\overline{\mathcal{T}}(\mathcal{B}) \neq \emptyset$.
(ii) $\underline{\mathcal{T}}(\mathcal{B})$ is a right bi-quasi ideal of $\mathcal{M}_1$, if $\underline{\mathcal{T}}(\mathcal{B}) \neq \emptyset$ and $\mathcal{T}$ is strong set-valued anti-homomorphism.

**Corollary 4.** Let $\mathcal{T}: \mathcal{M}_1 \to \mathcal{P}^*(\mathcal{M}_2)$ be a set-valued anti-homomorphism. If $\mathcal{B}$ is a bi-quasi ideal of $\mathcal{M}_2$. Then

(i) $\overline{\mathcal{T}}(\mathcal{B})$ is a bi-quasi ideal of $\mathcal{M}_1$, if $\overline{\mathcal{T}}(\mathcal{B}) \neq \emptyset$.
(ii) $\underline{\mathcal{T}}(\mathcal{B})$ is a bi-quasi ideal of $\mathcal{M}_1$, if $\underline{\mathcal{T}}(\mathcal{B}) \neq \emptyset$ and $\mathcal{T}$ is strong set-valued anti-homomorphism.

**Theorem 5.9.** Let $\mathcal{M}_1$ and $\mathcal{M}_2$ be two $\Gamma$-semigroups and $\mathcal{T}: \mathcal{M}_1 \to \mathcal{P}^*(\mathcal{M}_2)$ be a set-valued anti-homomorphism. If $Q$ is a left quasi-interior ideal of $\mathcal{M}_2$, then

(i) $\overline{\mathcal{T}}(Q)$ is a left quasi-interior ideal of $\mathcal{M}_1$, if $\mathcal{T}(x) \neq \emptyset, \forall\, x \in \mathcal{M}_1$ and $\overline{\mathcal{T}}(Q) \neq \emptyset$.
(ii) $\underline{\mathcal{T}}(Q)$ is a left quasi-interior ideal of $\mathcal{M}_1$, if $\underline{\mathcal{T}}(Q) \neq \emptyset$ and $\mathcal{T}$ is strong set-valued anti-homomorphism.

**Proof.** (i) Since $Q$ is a left quasi-interior ideal of $\mathcal{M}_2$, $\mathcal{M}_2\Gamma Q\Gamma\mathcal{M}_2\Gamma Q \subseteq Q$. By proposition $\overline{\mathcal{T}}(\mathcal{M}_2) = \mathcal{M}_1$. We have, $\mathcal{M}_1\Gamma\overline{\mathcal{T}}(Q)\Gamma\mathcal{M}_1\Gamma\overline{\mathcal{T}}(Q) = \overline{\mathcal{T}}(\mathcal{M}_2)\Gamma\overline{\mathcal{T}}(Q)\Gamma\overline{\mathcal{T}}(\mathcal{M}_2)\Gamma\overline{\mathcal{T}}(Q) \subseteq \overline{\mathcal{T}}(Q\Gamma\mathcal{M}_2\Gamma Q\Gamma\mathcal{M}_2) = \overline{\mathcal{T}}(\mathcal{M}_2\Gamma Q\Gamma\mathcal{M}_2\Gamma Q) \subseteq \overline{\mathcal{T}}(Q)$. Therefore $\overline{\mathcal{T}}(Q)$ is a left quasi-interior ideal of $\mathcal{M}_1$.

(ii) Given $Q$ is a left quasi-interior ideal of $\mathcal{M}_2$, $\mathcal{M}_2\Gamma Q\Gamma\mathcal{M}_2\Gamma Q \subseteq Q$. By proposition $\underline{\mathcal{T}}(\mathcal{M}_2) = \mathcal{M}_1$. Now, $\underline{\mathcal{T}}(Q)\Gamma\mathcal{M}_1\,\Gamma\underline{\mathcal{T}}(Q)\Gamma\mathcal{M}_1 = \underline{\mathcal{T}}(Q)\Gamma\underline{\mathcal{T}}(\mathcal{M}_2)\Gamma\underline{\mathcal{T}}(Q)\Gamma\underline{\mathcal{T}}(\mathcal{M}_2) \subseteq \underline{\mathcal{T}}(Q\Gamma\mathcal{M}_2\Gamma Q\Gamma\mathcal{M}_2) = \underline{\mathcal{T}}(\mathcal{M}_2\Gamma Q\Gamma\mathcal{M}_2\Gamma Q) \subseteq \underline{\mathcal{T}}(Q)$. Thus $\underline{\mathcal{T}}(\mathcal{M}_2\Gamma Q\Gamma\mathcal{M}_2\Gamma Q) \subseteq \underline{\mathcal{T}}(Q)$. Hence $\underline{\mathcal{T}}(Q)$ is a left quasi-interior ideal of $\mathcal{M}_1$.

**Corollary 5.** Let $\mathcal{T}: \mathcal{M}_1 \to \mathcal{P}^*(\mathcal{M}_2)$ be a set-valued anti-homomorphism. If $Q$ is a right quasi-interior ideal of $\mathcal{M}_2$. Then

(i) $\overline{\mathcal{T}}(Q)$ is a right quasi-interior ideal of $\mathcal{M}_1$, if $\overline{\mathcal{T}}(Q) \neq \emptyset$.
(ii) $\underline{\mathcal{T}}(Q)$ is a right quasi-interior ideal of $\mathcal{M}_1$, if $\underline{\mathcal{T}}(Q) \neq \emptyset$ and $\mathcal{T}$ is strong set-valued anti-homomorphism.

**Corollary 6.** Let $\mathcal{T}: \mathcal{M}_1 \to \mathcal{P}^*(\mathcal{M}_2)$ be a set-valued anti-homomorphism. If $\mathcal{Q}$ is a quasi-interior ideal of $\mathcal{M}_2$. Then

(i) $\overline{\mathcal{T}}(\mathcal{Q})$ is a quasi-interior ideal of $\mathcal{M}_1$, if $\overline{\mathcal{T}}(\mathcal{Q}) \neq \emptyset$.
(ii) $\underline{\mathcal{T}}(\mathcal{Q})$ is a quasi-interior ideal of $\mathcal{M}_1$, if $\underline{\mathcal{T}}(\mathcal{Q}) \neq \emptyset$ and $\mathcal{T}$ is strong set-valued anti-homomorphism.

**Theorem 5.10.** Let $\mathcal{M}_1$ and $\mathcal{M}_2$ be two $\Gamma$-semigroups and $\mathcal{T}: \mathcal{M}_1 \to \mathcal{P}^*(\mathcal{M}_2)$ be a set-valued anti-homomorphism. If $\mathcal{B}$ is a bi-quasi-interior ideal of $\mathcal{M}_2$, then

(i) $\overline{\mathcal{T}}(\mathcal{B})$ is a bi-quasi-interior ideal of $\mathcal{M}_1$, if $\mathcal{T}(x) \neq \emptyset, \forall\, x \in \mathcal{M}_1$ and $\overline{\mathcal{T}}(\mathcal{B}) \neq \emptyset$.
(ii) $\underline{\mathcal{T}}(\mathcal{B})$ is a bi-quasi-interior ideal of $\mathcal{M}_1$, if $\underline{\mathcal{T}}(\mathcal{B}) \neq \emptyset$ and $\mathcal{T}$ is strong set-valued anti-homomorphism.

**Proof.** (i) Since $\mathcal{B}$ is a bi-quasi-interior ideal of $\mathcal{M}_2$, $\mathcal{B}\Gamma\mathcal{M}_2\Gamma\mathcal{B}\Gamma\mathcal{M}_2\Gamma\mathcal{B} \subseteq \mathcal{B}$. By proposition $\overline{\mathcal{T}}(\mathcal{M}_2) = \mathcal{M}_1$. We have, $\overline{\mathcal{T}}(\mathcal{B})\Gamma\mathcal{M}_1\Gamma\overline{\mathcal{T}}(\mathcal{B})\Gamma\mathcal{M}_1\Gamma\overline{\mathcal{T}}(\mathcal{B}) =$
$\overline{\mathcal{T}}(\mathcal{B})\Gamma\overline{\mathcal{T}}(\mathcal{M}_2)\Gamma\overline{\mathcal{T}}(\mathcal{B})\Gamma\overline{\mathcal{T}}(\mathcal{M}_2)\Gamma\overline{\mathcal{T}}(\mathcal{B}) \subseteq \overline{\mathcal{T}}(\mathcal{B}\Gamma\mathcal{M}_2\Gamma\mathcal{B}\Gamma\mathcal{M}_2\Gamma\mathcal{B}) \subseteq \overline{\mathcal{T}}(\mathcal{B})$. Therefore $\overline{\mathcal{T}}(\mathcal{B})$ is a bi-quasi-interior ideal of $\mathcal{M}_1$.

(ii) Given $\mathcal{B}$ is a bi-quasi-interior ideal of $\mathcal{M}_2$, $\mathcal{B}\Gamma\mathcal{M}_2\Gamma\mathcal{B}\Gamma\mathcal{M}_2\Gamma\mathcal{B} \subseteq \mathcal{B}$. By proposition $\underline{\mathcal{T}}(\mathcal{M}_2) = \mathcal{M}_1$. Now, $\underline{\mathcal{T}}(\mathcal{B})\Gamma\mathcal{M}_1\Gamma\,\underline{\mathcal{T}}(\mathcal{B})\Gamma\mathcal{M}_1\Gamma\underline{\mathcal{T}}(\mathcal{B}) =$
$\underline{\mathcal{T}}(\mathcal{B})\Gamma\underline{\mathcal{T}}(\mathcal{M}_2)\Gamma\underline{\mathcal{T}}(\mathcal{B})\Gamma\underline{\mathcal{T}}(\mathcal{M}_2)\Gamma\underline{\mathcal{T}}(\mathcal{B}) \subseteq \underline{\mathcal{T}}(\mathcal{M}_2\Gamma\mathcal{B}\Gamma\mathcal{M}_2 \cap \mathcal{B}\Gamma\mathcal{M}_2\Gamma\mathcal{B}) \subseteq \underline{\mathcal{T}}(\mathcal{B})$. Thus $\underline{\mathcal{T}}(\mathcal{B}\Gamma\mathcal{M}_2\Gamma\mathcal{B}\Gamma\mathcal{M}_2\Gamma\mathcal{B}) \subseteq \underline{\mathcal{T}}(\mathcal{B})$. Hence $\underline{\mathcal{T}}(\mathcal{B})$ is a bi-quasi-interior ideal of $\mathcal{M}_1$.

**Theorem 5.11.** Consider $\mathcal{M}_1$ and $\mathcal{M}_2$ are two $\Gamma$-semigroups and define $\mathcal{T}: \mathcal{M}_1 \to \mathcal{P}^*(\mathcal{M}_2)$ is a set-valued anti-homomorphism. If $\mathcal{B}$ is a prime bi-ideal of $\mathcal{M}_2$ then

(i) $\overline{\mathcal{T}}(\mathcal{B})$ is a prime bi-ideal of $\mathcal{M}_1$, if $\overline{\mathcal{T}}(\mathcal{B}) \neq \emptyset$.
(ii) $\underline{\mathcal{T}}(\mathcal{B})$ is a prime bi-ideal of $\mathcal{M}_1$, if $\underline{\mathcal{T}}(\mathcal{B}) \neq \emptyset$ and $\mathcal{T}$ is strong set-valued anti-homomorphism.

**Proof.** (i) By theorem 5.4.(i), $\overline{\mathcal{T}}(\mathcal{B})$ is a bi-ideal of $\mathcal{M}_1$. Assume $a\alpha b \in \overline{\mathcal{T}}(\mathcal{B})$ and $\mathcal{T}(a\alpha b) \in \mathcal{B}$. Then $\mathcal{T}(a\alpha b) \cap \mathcal{B} \neq \emptyset$. We have, $\mathcal{T}(a\alpha b) \supseteq \mathcal{T}(b)\alpha\mathcal{T}(a), \forall\, a,b \in \mathcal{M}_1$ and $\alpha \in \Gamma$. There exists $u \in \mathcal{T}(a), v \in \mathcal{T}(b)$ such that $v\alpha u \subseteq \mathcal{T}(a\alpha b)$. (i.e) $v\alpha u \in \mathcal{B}$. Since $\mathcal{B}$ is a prime bi-ideal of $\mathcal{M}_2$, $u \in \mathcal{B}$ or $v \in \mathcal{B}$. Thus $a \in \overline{\mathcal{T}}(\mathcal{B})$ or $b \in \overline{\mathcal{T}}(\mathcal{B})$. Hence $\overline{\mathcal{T}}(\mathcal{B})$ is a prime bi-ideal of $\mathcal{M}_1$.

(ii) By theorem 5.4.(ii), $\underline{\mathcal{T}}(\mathcal{B})$ is a bi-ideal of $\mathcal{M}_1$. Suppose $\underline{\mathcal{T}}(\mathcal{B})$ is not prime bi-ideal of $\mathcal{M}_1$. Then there exists $a, b \in \mathcal{M}_1$ such that $a\alpha b \in \underline{\mathcal{T}}(\mathcal{B})$ but $a \notin \underline{\mathcal{T}}(\mathcal{B})$, $b \notin \underline{\mathcal{T}}(\mathcal{B})$. There exists $u \in \mathcal{T}(a), v \in \mathcal{T}(b)$ but $u, v \notin \mathcal{B}$. Thus $v\alpha u \in \mathcal{T}(b)\alpha\mathcal{T}(a) \subseteq \mathcal{T}(a\alpha b) \subseteq \mathcal{B}, \forall\, a, b \in \mathcal{M}_1$ and $\alpha \in \Gamma$. Since $\mathcal{B}$ is a prime bi-ideal of $\mathcal{M}_2$, $u \in \mathcal{B}$ or $v \in \mathcal{B}$, which is a contradiction. Therefore $\underline{\mathcal{T}}(\mathcal{B})$ is a prime bi-ideal of $\mathcal{M}_1$.

**Example 13.** Consider the $\Gamma$-semigroup $\mathcal{M}$ of Example 4. Define a set-valued mapping $\mathcal{T}: \mathcal{M} \to \mathcal{P}^*(\mathcal{M})$ by $\mathcal{T}(x_1) = \{x_1, x_2, x_3, x_4\}, \mathcal{T}(x_2) = \{x_1, x_3\}, \mathcal{T}(x_3) = \{x_3\}, \mathcal{T}(x_4) = \{x_4\}$. Then $\mathcal{T}$ is a set-valued anti-homomorphism. Assume the subsets $\{x_1\}$, $\{x_2\}, \{x_3, x_4\}$. Let $\mathcal{B} = \{x_1, x_2, x_3\} \subset \mathcal{M}$. Then $\mathcal{B}$ is a bi-ideal of $\mathcal{M}$. $\overline{\mathcal{T}}(\mathcal{B}) = \{x_1, x_2, x_3, x_4\}$, $\mathcal{B}\Gamma\mathcal{M}\Gamma\mathcal{B} = \{x_1, x_3\}$, $\overline{\mathcal{T}}(\mathcal{B}\Gamma\mathcal{M}\Gamma\mathcal{B}) = \{x_1, x_3, x_4\}$. Hence $\overline{\mathcal{T}}(\mathcal{B}\Gamma\mathcal{M}\Gamma\mathcal{B}) \subseteq$

$\overline{\mathcal{T}}(\mathcal{B})$. Therefore $\overline{\mathcal{T}}(\mathcal{B})$ is bi-ideal of $\mathcal{M}$. It is obvious that $\mathcal{B}$ is a prime bi-ideal of $\Gamma$-semigroup $\mathcal{M}$. Here, the set $\overline{\mathcal{T}}(\mathcal{B})$ is a prime bi-ideal of $\mathcal{M}$ for $x_1 \alpha x_3 \in \overline{\mathcal{T}}(\mathcal{B})$ as $x_1 \in \overline{\mathcal{T}}(\mathcal{B})$, $x_3 \in \overline{\mathcal{T}}(\mathcal{B})$.

**Example 14.** Consider the $\Gamma$-semigroups $\mathcal{M}_1$ and $\mathcal{M}_2$ of Example 5. Define a set-valued mapping $\mathcal{T} \colon \mathcal{M}_1 \to \mathcal{P}^*(\mathcal{M}_2)$ Assume $\mathcal{T}(1) = \{c\}$, $\mathcal{T}(2) = \{a\}$. Then $\mathcal{T}$ is a strong set-valued anti-homomorphism. Assuming the subsets $\{a\}, \{c\}, \{a, b\}$. Let $\mathcal{B} = \{b, c\} \subset \mathcal{M}_2$. Then $\mathcal{B}$ is a bi-ideal of $\mathcal{M}_2$. $\underline{\mathcal{T}}(\mathcal{B}) = \{b, c\}$, $\mathcal{B}\Gamma\mathcal{M}_2\Gamma\mathcal{B} = \{b\}$ and $\underline{\mathcal{T}}(\mathcal{B}\Gamma\mathcal{M}_2\Gamma\mathcal{B}) = \{b\}$. Hence $\underline{\mathcal{T}}(\mathcal{B}\Gamma\mathcal{M}_2\Gamma\mathcal{B}) \subseteq \underline{\mathcal{T}}(\mathcal{B})$. Therefore $\underline{\mathcal{T}}(\mathcal{B})$ is a bi-ideal of $\mathcal{M}_1$. It is clear that $\mathcal{B}$ is a prime bi-ideal of $\Gamma$-semigroup $\mathcal{M}_1$. Also, the set $\underline{\mathcal{T}}(\mathcal{B})$ is a prime bi-ideal of $\mathcal{M}_1$ for $a\alpha c \in \underline{\mathcal{T}}(\mathcal{B})$ as $c \in \underline{\mathcal{T}}(\mathcal{B})$, $a \notin \underline{\mathcal{T}}(\mathcal{B})$.

**Theorem 5.12.** Define $\mathcal{T} \colon \mathcal{M}_1 \to \mathcal{P}^*(\mathcal{M}_2)$ a set-valued anti-homomorphism where $\mathcal{M}_1$ and $\mathcal{M}_2$ are two $\Gamma$-semigroups. Suppose $\mathcal{J}$ is a prime interior ideal of $\mathcal{M}_2$. Then

   (i) $\overline{\mathcal{T}}(\mathcal{J})$ is a prime interior ideal of $\mathcal{M}_1$, if $\overline{\mathcal{T}}(\mathcal{J}) \neq \emptyset$.
   (ii) $\underline{\mathcal{T}}(\mathcal{J})$ is a prime interior ideal of $\mathcal{M}_1$, if $\underline{\mathcal{T}}(\mathcal{J}) \neq \emptyset$ and $\mathcal{T}$ is strong set-valued anti-homomorphism.

**Proof.** (i) By theorem 5.5.(i), $\overline{\mathcal{T}}(\mathcal{J})$ is an interior ideal of $\mathcal{M}_1$. Let $x\alpha y \in \overline{\mathcal{T}}(\mathcal{J})$. Then $\mathcal{T}(x\alpha y) \cap \mathcal{J} \neq \emptyset$. We have $\mathcal{T}(x\alpha y) \supseteq \mathcal{T}(y)\alpha\mathcal{T}(x), \forall\, x, y \in \mathcal{M}_1$ and $\alpha \in \Gamma$. There exists $u \in \mathcal{T}(x), v \in \mathcal{T}(y)$ such that $v\alpha u \subseteq \mathcal{T}(x\alpha y)$. So, $v\alpha u \in \mathcal{J}$. Since $\mathcal{J}$ is a prime interior ideal of $\mathcal{M}_2$, $u \in \mathcal{J}$ or $v \in \mathcal{J}$. Therefore $x \in \overline{\mathcal{T}}(\mathcal{J})$ or $y \in \overline{\mathcal{T}}(\mathcal{J})$. Thus $\overline{\mathcal{T}}(\mathcal{J})$ is a prime interior ideal of $\mathcal{M}_1$.

(ii) By theorem 5.5.(ii), $\underline{\mathcal{T}}(\mathcal{J})$ is an interior ideal of $\mathcal{M}_1$. Suppose $\underline{\mathcal{T}}(\mathcal{J})$ is not prime interior ideal of $\mathcal{M}_1$. There exists $x, y \in \mathcal{M}_1$ such that $x\alpha y \in \underline{\mathcal{T}}(\mathcal{J})$ but $x \notin \underline{\mathcal{T}}(\mathcal{J})$, $y \notin \underline{\mathcal{T}}(\mathcal{J})$. Then $u \in \mathcal{T}(x), v \in \mathcal{T}(y)$ but $u, v \notin \mathcal{J}$. Thus $v\alpha u \in \mathcal{T}(y)\alpha\mathcal{T}(x) \subseteq \mathcal{T}(x\alpha y) \subseteq \mathcal{J}, \forall\, x, y \in \mathcal{M}_1$ and $\alpha \in \Gamma$. Since $\mathcal{J}$ is a prime interior ideal of $\mathcal{M}_2$, $u \in \mathcal{J}$ or $v \in \mathcal{J}$ which contradicts itself. Hence $\underline{\mathcal{T}}(\mathcal{J})$ is a prime interior ideal of $\mathcal{M}_1$.

**Example 15.** Consider the $\Gamma$-semigroup $\mathcal{M}$ of Example 4. Define a set-valued mapping $\mathcal{T} \colon \mathcal{M} \to \mathcal{P}^*(\mathcal{M})$ by $\mathcal{T}(x_1) = \{x_1, x_2, x_3, x_4\}, \mathcal{T}(x_2) = \{x_1, x_3\}, \mathcal{T}(x_3) = \{x_3\}, \mathcal{T}(x_4) = \{x_4\}$. Then $\mathcal{T}$ is a set-valued anti-homomorphism. Assume the subsets $\{x_1, x_2, x_3\}, \{x_4\}$. Let $\mathcal{J} = \{x_1, x_3, x_4\} \subset \mathcal{M}$. Then $\mathcal{J}$ is an interior ideal of $\mathcal{M}$. $\overline{\mathcal{T}}(\mathcal{J}) = \{x_1, x_2, x_3, x_4\}$, $\mathcal{M}\Gamma\mathcal{J}\Gamma\mathcal{M} = \{x_1, x_3\}$, $\overline{\mathcal{T}}(\mathcal{M}\Gamma\mathcal{J}\Gamma\mathcal{M}) = \{x_1, x_2, x_3\}$. Hence $\overline{\mathcal{T}}(\mathcal{M}\Gamma\mathcal{J}\Gamma\mathcal{M}) \subseteq \overline{\mathcal{T}}(\mathcal{J})$. Therefore $\overline{\mathcal{T}}(\mathcal{J})$ is an interior ideal of $\mathcal{M}$. It is obvious that $\mathcal{J}$ is a prime interior ideal of $\Gamma$-semigroup $\mathcal{M}$. Here, the set $\overline{\mathcal{T}}(\mathcal{J})$ is a prime interior ideal of $\mathcal{M}$ for $x_1\alpha x_4 \in \overline{\mathcal{T}}(\mathcal{J})$ as $x_1 \in \overline{\mathcal{T}}(\mathcal{J})$, $x_4 \in \overline{\mathcal{T}}(\mathcal{J})$.

**Example 16.** Consider the $\Gamma$-semigroups $\mathcal{M}_1$ and $\mathcal{M}_2$ of Example 5. Define a set-valued mapping $\mathcal{T} \colon \mathcal{M}_1 \to \mathcal{P}^*(\mathcal{M}_2)$ Assume $\mathcal{T}(1) = \{c\}$, $\mathcal{T}(2) = \{a\}$. Then $\mathcal{T}$ is a strong set-valued anti-homomorphism. Assuming the subsets $\{b\}, \{c\}, \{a, b\}$. Let $\mathcal{J} = \{b, c\} \subset \mathcal{M}_2$. Then $\mathcal{J}$ is an interior ideal of $\mathcal{M}_2$. $\underline{\mathcal{T}}(\mathcal{J}) = \{b, c\}$, $\mathcal{M}_2\Gamma\mathcal{J}\Gamma\mathcal{M}_2 = \{b, c\}$ and $\underline{\mathcal{T}}(\mathcal{M}_2\Gamma\mathcal{J}\Gamma\mathcal{M}_2) = \{b, c\}$. Hence $\underline{\mathcal{T}}(\mathcal{M}_2\Gamma\mathcal{J}\Gamma\mathcal{M}_2) \subseteq \underline{\mathcal{T}}(\mathcal{J})$. Therefore $\underline{\mathcal{T}}(\mathcal{J})$ is an interior ideal of $\mathcal{M}_1$. It is clear that $\mathcal{J}$ is a prime interior ideal of $\Gamma$-semigroup $\mathcal{M}_1$. Also, the set $\underline{\mathcal{T}}(\mathcal{J})$ is a prime interior ideal of $\mathcal{M}_1$ for $b\alpha c \in \underline{\mathcal{T}}(\mathcal{J})$ as $b \in \underline{\mathcal{T}}(\mathcal{J})$ or $c \in \underline{\mathcal{T}}(\mathcal{J})$.

**Theorem 5.13.** Consider two Γ-semigroups $\mathcal{M}_1$ and $\mathcal{M}_2$. Defining $\mathcal{T}: \mathcal{M}_1 \to \mathcal{P}^*(\mathcal{M}_2)$ a set-valued anti-homomorphism. Let $Q$ be a prime quasi-ideal of $\mathcal{M}_2$ then

(i)  $\overline{\mathcal{T}}(Q)$ is a prime quasi-ideal of $\mathcal{M}_1$, if $\overline{\mathcal{T}}(Q) \neq \emptyset$.
(ii) $\underline{\mathcal{T}}(Q)$ is a prime quasi-ideal of $\mathcal{M}_1$, if $\underline{\mathcal{T}}(Q) \neq \emptyset$ and $\mathcal{T}$ is strong set-valued anti-homomorphism.

**Proof.** (i) By theorem 5.6.(i), $\overline{\mathcal{T}}(Q)$ is a quasi-ideal of $\mathcal{M}_1$. Take aαb ∈ $\overline{\mathcal{T}}(Q)$ and $\mathcal{T}$(aαb) ∈ $Q$. Then $\mathcal{T}$(aαb) ∩ $Q \neq \emptyset$. We have, $\mathcal{T}$(aαb) ⊇ $\mathcal{T}$(b)α$\mathcal{T}$(a), ∀ a, b ∈ $\mathcal{M}_1$ and α ∈ Γ. There exists u ∈ $\mathcal{T}$(a), v ∈ $\mathcal{T}$(b) such that vαu ⊆ $\mathcal{T}$(aαb). (i.e) vαu ∈ $Q$. Since $Q$ is prime quasi-ideal of $\mathcal{M}_2$, u ∈ $Q$ or v ∈ $Q$. Therefore a ∈ $\overline{\mathcal{T}}(Q)$ or b ∈ $\overline{\mathcal{T}}(Q)$. Hence $\overline{\mathcal{T}}(Q)$ is a prime quasi-ideal of $\mathcal{M}_1$.

(ii) By theorem 5.6.(ii), $\underline{\mathcal{T}}(Q)$ is a quasi-ideal of $\mathcal{M}_1$. Suppose $\underline{\mathcal{T}}(Q)$ is not prime quasi-ideal of $\mathcal{M}_1$. There exists a, b ∈ $\mathcal{M}_1$ such that aαb ∈ $\underline{\mathcal{T}}(Q)$ but a ∉ $\underline{\mathcal{T}}(Q)$, b ∉ $\underline{\mathcal{T}}(Q)$. Then, as u ∈ $\mathcal{T}$(a), v ∈ $\mathcal{T}$(b), it follows that u, v ∉ $Q$. Thus vαu ∈ $\mathcal{T}$(b)α$\mathcal{T}$(a) ⊆ $\mathcal{T}$(aαb) ⊆ $Q$, ∀ a, b ∈ $\mathcal{M}_1$ and α ∈ Γ. Since $Q$ is a prime quasi-ideal of $\mathcal{M}_2$, u ∈ $Q$ or v ∈ $Q$ which is incongruous. Hence $\underline{\mathcal{T}}(Q)$ is a prime quasi-ideal of $\mathcal{M}_1$.

**Example 17.** Consider the Γ-semigroup $\mathcal{M}$ of Example 4. Define a set-valued mapping $\mathcal{T}: \mathcal{M} \to \mathcal{P}^*(\mathcal{M})$ by $\mathcal{T}(x_1) = \{x_1, x_2, x_3, x_4\}, \mathcal{T}(x_2) = \{x_1, x_3\}, \mathcal{T}(x_3) = \{x_3\}, \mathcal{T}(x_4) = \{x_4\}$. Then $\mathcal{T}$ is a set-valued anti-homomorphism. Assume the subsets $\{x_1, x_2\}, \{x_3\}, \{x_4\}$. Let $Q = \{x_1, x_3, x_4\} \subset \mathcal{M}$. Then $Q$ is a quasi-ideal of $\mathcal{M}$. $\overline{\mathcal{T}}(Q) = \{x_1, x_2, x_3, x_4\}$, $Q\Gamma\mathcal{M} \cap \mathcal{M}\Gamma Q = \{x_1, x_3\}$, $\overline{\mathcal{T}}(Q\Gamma\mathcal{M} \cap \mathcal{M}\Gamma Q) = \{x_1, x_2, x_3\}$. Hence $\overline{\mathcal{T}}(Q\Gamma\mathcal{M} \cap \mathcal{M}\Gamma Q) \subseteq \overline{\mathcal{T}}(Q)$. Therefore $\overline{\mathcal{T}}(Q)$ is a quasi-ideal of $\mathcal{M}$. It is obvious that $Q$ is a prime quasi-ideal of Γ-semigroup $\mathcal{M}$. Here, the set $\overline{\mathcal{T}}(Q)$ is a prime quasi-ideal of $\mathcal{M}$ for $x_2 \alpha x_3 \in \overline{\mathcal{T}}(Q)$ as $x_2 \in \overline{\mathcal{T}}(Q), x_3 \in \overline{\mathcal{T}}(Q)$.

**Example 18.** Consider the Γ-semigroups $\mathcal{M}_1$ and $\mathcal{M}_2$ of Example 6. Define a set-valued mapping $\mathcal{T}: \mathcal{M}_1 \to \mathcal{P}^*(\mathcal{M}_2)$ by $\mathcal{T}(x) = \mathcal{T}(z) = \{c\}, \mathcal{T}(y) = \{b, c\}$. Then $\mathcal{T}$ is a strong set-valued anti-homomorphism. Assuming the subsets $\{c\}, \{b, c\}$ and Let $Q = \{a, b, c\} \subset \mathcal{M}_2$. Then $Q$ is a quasi-ideal of $\mathcal{M}_2$. $\underline{\mathcal{T}}(Q) = \{b, c\}$, $Q\Gamma\mathcal{M}_2 \cap \mathcal{M}_2\Gamma Q = \{a, c\}$ and $\underline{\mathcal{T}}(Q\Gamma\mathcal{M}_2 \cap \mathcal{M}_2\Gamma Q) = \{c\}$. Hence $\underline{\mathcal{T}}(Q\Gamma\mathcal{M}_2 \cap \mathcal{M}_2\Gamma Q) \subseteq \underline{\mathcal{T}}(Q)$. Therefore $\underline{\mathcal{T}}(Q)$ is a quasi-ideal of $\mathcal{M}_1$. It is obvious that $Q$ is a prime quasi-ideal of Γ-semigroup $\mathcal{M}_1$. Also, the set $\underline{\mathcal{T}}(Q)$ is not a prime quasi-ideal of $\mathcal{M}_1$ for $a\alpha b \notin \underline{\mathcal{T}}(Q)$ as a ∉ $\underline{\mathcal{T}}(Q)$ or b ∈ $\underline{\mathcal{T}}(Q)$.

**Theorem 5.14.** Consider $\mathcal{M}_1$ and $\mathcal{M}_2$ are two Γ-semigroups and define $\mathcal{T}: \mathcal{M}_1 \to \mathcal{P}^*(\mathcal{M}_2)$ is a set-valued anti-homomorphism. If $\mathcal{B}$ is a prime bi-interior ideal of $\mathcal{M}_2$ then

(i)  $\overline{\mathcal{T}}(\mathcal{B})$ is a prime bi-interior ideal of $\mathcal{M}_1$, if $\overline{\mathcal{T}}(\mathcal{B}) \neq \emptyset$.
(ii) $\underline{\mathcal{T}}(\mathcal{B})$ is a prime bi-interior ideal of $\mathcal{M}_1$, if $\underline{\mathcal{T}}(\mathcal{B}) \neq \emptyset$ and $\mathcal{T}$ is strong set-valued anti-homomorphism.

**Proof.** (i) By theorem 5.7.(i), $\overline{\mathcal{T}}(\mathcal{B})$ is a bi-interior ideal of $\mathcal{M}_1$. Assume aαb ∈ $\overline{\mathcal{T}}(\mathcal{B})$ and $\mathcal{T}$(aαb) ∈ $\mathcal{B}$. Then $\mathcal{T}$(aαb) ∩ $\mathcal{B} \neq \emptyset$. We have, $\mathcal{T}$(aαb) ⊇ $\mathcal{T}$(b)α$\mathcal{T}$(a), ∀ a, b ∈ $\mathcal{M}_1$ and α ∈ Γ. There exists u ∈ $\mathcal{T}$(a), v ∈ $\mathcal{T}$(b) such that vαu ⊆ $\mathcal{T}$(aαb). (i.e) vαu ∈ $\mathcal{B}$. Since $\mathcal{B}$ is a prime bi-interior ideal of $\mathcal{M}_2$, $u \in \mathcal{B}$ or $v \in \mathcal{B}$. Thus $a \in \overline{\mathcal{T}}(\mathcal{B})$ or $b \in \overline{\mathcal{T}}(\mathcal{B})$. Hence $\overline{\mathcal{T}}(\mathcal{B})$ is a prime bi-interior ideal of $\mathcal{M}_1$.

(ii) By theorem 5.7.(ii), $\underline{\mathcal{T}}(\mathcal{B})$ is a bi-interior ideal of $\mathcal{M}_1$. Suppose $\underline{\mathcal{T}}(\mathcal{B})$ is not prime bi-interior ideal of $\mathcal{M}_1$. Then there exists $a, b \in \mathcal{M}_1$ such that $a\alpha b \in \underline{\mathcal{T}}(\mathcal{B})$ but $a \notin \underline{\mathcal{T}}(\mathcal{B})$, $b \notin \underline{\mathcal{T}}(\mathcal{B})$. There exists $u \in \mathcal{T}(a), v \in \mathcal{T}(b)$ but $u, v \notin \mathcal{B}$. Thus $v\alpha u \in \mathcal{T}(b)\alpha\mathcal{T}(a) \subseteq \mathcal{T}(a\alpha b) \subseteq \mathcal{B}, \forall a, b \in \mathcal{M}_1$ and $\alpha \in \Gamma$. Since $\mathcal{B}$ is a prime bi-interior ideal of $\mathcal{M}_2$, $u \in \mathcal{B}$ or $v \in \mathcal{B}$, which is a contradiction. Therefore, $\underline{\mathcal{T}}(\mathcal{B})$ is a prime bi-interior ideal of $\mathcal{M}_1$.

**Example 19.** Consider the $\Gamma$-semigroup $\mathcal{M}$ of Example 4. Define a set-valued mapping $\mathcal{T}: \mathcal{M} \to \mathcal{P}^*(\mathcal{M})$ by $\mathcal{T}(x_1) = \{x_1, x_2, x_3, x_4\}, \mathcal{T}(x_2) = \{x_1, x_3\}, \mathcal{T}(x_3) = \{x_3\}, \mathcal{T}(x_4) = \{x_4\}$. Then $\mathcal{T}$ is a set-valued anti-homomorphism. $\mathcal{M} = \{x_1, x_2, x_3, x_4\}$ and $\Gamma = \{\alpha\}$ and taking the subsets $\{x_1, x_2\}, \{x_1, x_3\}, \{x_3, x_4\}$. Then for $\mathcal{B} = \{x_1, x_3\} \subseteq \mathcal{M}$, Then $\mathcal{B}$ is a bi-interior ideal of $\mathcal{M}$. $\overline{\mathcal{T}}(\mathcal{B}) = \{x_1, x_2, x_3, x_4\}$, $\underline{\mathcal{T}}(\mathcal{B}) = \{x_1, x_3\}$, $\mathcal{M}\Gamma\mathcal{B}\Gamma\mathcal{M} \cap \mathcal{B}\Gamma\mathcal{M}\Gamma\mathcal{B} = \{x_1, x_3\}$, $\overline{\mathcal{T}}(\mathcal{M}\Gamma\mathcal{B}\Gamma\mathcal{M} \cap \mathcal{B}\Gamma\mathcal{M}\Gamma\mathcal{B}) = \{x_1, x_2, x_3, x_4\}$. Hence $\overline{\mathcal{T}}(\mathcal{M}\Gamma\mathcal{B}\Gamma\mathcal{M} \cap \mathcal{B}\Gamma\mathcal{M}\Gamma\mathcal{B}) \subseteq \overline{\mathcal{T}}(\mathcal{B})$. Therefore $\overline{\mathcal{T}}(\mathcal{B})$ is a bi-interior ideal of $\mathcal{M}$. It is obvious that $\mathcal{B}$ is a prime bi-interior ideal of $\Gamma$-semigroup $\mathcal{M}$. Here, the set $\overline{\mathcal{T}}(\mathcal{B})$ is a prime bi-interior ideal of $\mathcal{M}$ for $x_2 \alpha x_4 \in \overline{\mathcal{T}}(\mathcal{B})$ as $x_2 \in \overline{\mathcal{T}}(\mathcal{B}), x_4 \in \overline{\mathcal{T}}(\mathcal{B})$.

**Example 20.** Consider the $\Gamma$-semigroups $\mathcal{M}_1$ and $\mathcal{M}_2$ of Example 6. Define a set-valued mapping $\mathcal{T}: \mathcal{M}_1 \to \mathcal{P}^*(\mathcal{M}_2)$ by $\mathcal{T}(x) = \mathcal{T}(z) = \{c\}, \mathcal{T}(y) = \{b, c\}$. Then $\mathcal{T}$ is a strong set-valued anti-homomorphism. Assuming the subsets $\{a,b\}\{b\},\{a,c\}$ and Let $\mathcal{B} = \{a, b, c\} \subset \mathcal{M}_2$. Then $\mathcal{B}$ is a bi-interior ideal of $\mathcal{M}_2$. $\underline{\mathcal{T}}(\mathcal{B}) = \{a, b, c\}$, $\mathcal{M}\Gamma\mathcal{B}\Gamma\mathcal{M} \cap \mathcal{B}\Gamma\mathcal{M}\Gamma\mathcal{B} = \{a, c\}$ and $\underline{\mathcal{T}}(\mathcal{M}\Gamma\mathcal{B}\Gamma\mathcal{M} \cap \mathcal{B}\Gamma\mathcal{M}\Gamma\mathcal{B}) = \{a, c\}$. Hence $\underline{\mathcal{T}}(\mathcal{M}\Gamma\mathcal{B}\Gamma\mathcal{M} \cap \mathcal{B}\Gamma\mathcal{M}\Gamma\mathcal{B}) \subseteq \underline{\mathcal{T}}(\mathcal{B})$. Therefore $\underline{\mathcal{T}}(\mathcal{B})$ is a bi- interior ideal of $\mathcal{M}_1$. It is obvious that $\mathcal{B}$ is a prime bi-interior ideal of $\Gamma$-semigroup $\mathcal{M}_1$. Also, the set $\underline{\mathcal{T}}(\mathcal{B})$ is a prime bi-interior ideal of $\mathcal{M}_1$ for $a\alpha c \in \underline{\mathcal{T}}(\mathcal{B})$ as $a \in \underline{\mathcal{T}}(\mathcal{B})$ or $b \in \underline{\mathcal{T}}(\mathcal{B})$.

**Theorem 5.15.** Define $\mathcal{T}: \mathcal{M}_1 \to \mathcal{P}^*(\mathcal{M}_2)$ a set-valued anti-homomorphism where $\mathcal{M}_1$ and $\mathcal{M}_2$ are two $\Gamma$-semigroups. $\mathcal{M}\Gamma\mathcal{B}\Gamma\mathcal{M} \cap \mathcal{B}\Gamma\mathcal{M}\Gamma\mathcal{B}$ Then

(i) $\overline{\mathcal{T}}(\mathcal{B})$ is a prime(left) bi-quasi ideal of $\mathcal{M}_1$, if $\overline{\mathcal{T}}(\mathcal{B}) \neq \emptyset$.
(ii) $\underline{\mathcal{T}}(\mathcal{B})$ is a prime(left) bi-quasi ideal of $\mathcal{M}_1$, if $\underline{\mathcal{T}}(\mathcal{B}) \neq \emptyset$ and $\mathcal{T}$ is strong set-valued anti-homomorphism.

**Proof.** (i) By theorem 5.8.(i), $\overline{\mathcal{T}}(\mathcal{B})$ is a (left)bi-quasi ideal of $\mathcal{M}_1$. Let $x\alpha y \in \overline{\mathcal{T}}(\mathcal{B})$. Then $\mathcal{T}(x\alpha y) \cap \mathcal{B} \neq \emptyset$. We have $\mathcal{T}(x\alpha y) \supseteq \mathcal{T}(y)\alpha\mathcal{T}(x), \forall x, y \in \mathcal{M}_1$ and $\alpha \in \Gamma$. There exists $u \in \mathcal{T}(x), v \in \mathcal{T}(y)$ such that $v\alpha u \subseteq \mathcal{T}(x\alpha y)$. So, $v\alpha u \in \mathcal{B}$. Since $\mathcal{B}$ is a prime(left) bi-quasi ideal of $\mathcal{M}_2$, $u \in \mathcal{B}$ or $v \in \mathcal{B}$. Therefore $x \in \overline{\mathcal{T}}(\mathcal{B})$ or $y \in \overline{\mathcal{T}}(\mathcal{B})$. Thus $\overline{\mathcal{T}}(\mathcal{B})$ is a prime(left) bi-quasi ideal of $\mathcal{M}_1$.

(ii) By theorem 5.8.(ii), $\underline{\mathcal{T}}(\mathcal{B})$ is a (left)bi-quasi ideal of $\mathcal{M}_1$. Suppose $\underline{\mathcal{T}}(\mathcal{B})$ is not prime (left)bi-quasi ideal of $\mathcal{M}_1$. There exists $x, y \in \mathcal{M}_1$ such that $x\alpha y \in \underline{\mathcal{T}}(\mathcal{B})$ but $x \notin \underline{\mathcal{T}}(\mathcal{B})$, $y \notin \underline{\mathcal{T}}(\mathcal{B})$. Then $u \in \mathcal{T}(x), v \in \mathcal{T}(y)$ but $u, v \notin \mathcal{J}$. Thus $v\alpha u \in \mathcal{T}(y)\alpha\mathcal{T}(x) \subseteq \mathcal{T}(x\alpha y) \subseteq \mathcal{B}, \forall x, y \in \mathcal{M}_1$ and $\alpha \in \Gamma$. Since $\mathcal{B}$ is a prime(left) bi-quasi ideal of $\mathcal{M}_2$, $u \in \mathcal{B}$ or $v \in \mathcal{B}$ which contradicts itself. Hence $\underline{\mathcal{T}}(\mathcal{B})$ is a prime(left) bi-quasi ideal of $\mathcal{M}_1$.

**Example 21.** Consider the Γ-semigroup $\mathcal{M}$ of Example 4. Define a set-valued mapping $\mathcal{T}: \mathcal{M} \to \mathcal{P}^*(\mathcal{M})$ by $\mathcal{T}(x_1) = \{x_1, x_2, x_3, x_4\}, \mathcal{T}(x_2) = \{x_1, x_3\}, \mathcal{T}(x_3) = \{x_3\}, \mathcal{T}(x_4) = \{x_4\}$. Then $\mathcal{T}$ is a set-valued anti-homomorphism. $\mathcal{M} = \{x_1, x_2, x_3, x_4\}$ and $\Gamma = \{\alpha\}$ and considering the subsets $\{x_1, x_2, x_3\}, \{x_3\}, \{x_1, x_4\}$. Let $\mathcal{B} = \{x_1, x_3\} \subseteq \mathcal{M}$. Then $\mathcal{B}$ is a (left)bi-quasi ideal of $\mathcal{M}$. $\overline{\mathcal{T}}(\mathcal{B}) = \{x_1, x_2, x_3, x_4\}$, $\mathcal{M}\Gamma\mathcal{B} \cap \mathcal{B}\Gamma\mathcal{M}\Gamma\mathcal{B} = \{x_1, x_3\}$, $\overline{\mathcal{T}}(\mathcal{M}\Gamma\mathcal{B} \cap \mathcal{B}\Gamma\mathcal{M}\Gamma\mathcal{B}) = \{x_1, x_2, x_3, x_4\}$. Hence $\overline{\mathcal{T}}(\mathcal{M}\Gamma\mathcal{B} \cap \mathcal{B}\Gamma\mathcal{M}\Gamma\mathcal{B}) \subseteq \overline{\mathcal{T}}(\mathcal{B})$. Therefore $\overline{\mathcal{T}}(\mathcal{B})$ is a (left)bi-quasi ideal of $\mathcal{M}$. It is obvious that $\mathcal{B}$ is a prime(left) bi-quasi ideal of Γ-semigroup $\mathcal{M}$. Here, the set $\overline{\mathcal{T}}(\mathcal{B})$ is a prime(left) bi-quasi ideal of $\mathcal{M}$ for $x_1 \alpha x_4 \in \overline{\mathcal{T}}(\mathcal{B})$ as $x_1 \in \overline{\mathcal{T}}(\mathcal{B}), x_4 \in \overline{\mathcal{T}}(\mathcal{B})$.

**Example 22.** Consider the Γ-semigroups $\mathcal{M}_1$ and $\mathcal{M}_2$ of Example 6. Define a set-valued mapping $\mathcal{T}: \mathcal{M}_1 \to \mathcal{P}^*(\mathcal{M}_2)$ by $\mathcal{T}(x) = \mathcal{T}(z) = \{c\}, \mathcal{T}(y) = \{b, c\}$. Then $\mathcal{T}$ is a strong set-valued anti-homomorphism. Assuming the subsets $\{c\}\{a\},\{b,c\}$ and Let $\mathcal{B} = \{a, c\} \subset \mathcal{M}_2$. Then $\mathcal{B}$ is a bi-quasi ideal of $\mathcal{M}_2$. Also, $\underline{\mathcal{T}}(\mathcal{B}) = \{a, c\}, \mathcal{M}\Gamma\mathcal{B} \cap \mathcal{B}\Gamma\mathcal{M}\Gamma\mathcal{B} = \{a, c\}$ and $\underline{\mathcal{T}}(\mathcal{M}\Gamma\mathcal{B} \cap \mathcal{B}\Gamma\mathcal{M}\Gamma\mathcal{B}) = \{a, c\}$. Hence $\underline{\mathcal{T}}(\mathcal{M}\Gamma\mathcal{B} \cap \mathcal{B}\Gamma\mathcal{M}\Gamma\mathcal{B}) \subseteq \underline{\mathcal{T}}(\mathcal{B})$. Therefore $\underline{\mathcal{T}}(\mathcal{B})$ is a bi-quasi ideal of $\mathcal{M}_1$. It is obvious that $\mathcal{B}$ is a prime bi-quasi ideal of Γ-semigroup $\mathcal{M}_1$. Also, the set $\underline{\mathcal{T}}(\mathcal{B})$ is a prime bi-quasi ideal of $\mathcal{M}_1$ for bαc $\in \underline{\mathcal{T}}(\mathcal{B})$ as b $\notin \underline{\mathcal{T}}(\mathcal{B})$ or c $\in \underline{\mathcal{T}}(\mathcal{B})$.

**Theorem 5.16.** Consider two Γ-semigroups $\mathcal{M}_1$ and $\mathcal{M}_2$. Defining $\mathcal{T}: \mathcal{M}_1 \to \mathcal{P}^*(\mathcal{M}_2)$ a set-valued anti-homomorphism. Let $Q$ be a prime(left) quasi-interior ideal of $\mathcal{M}_2$ then

(i) $\overline{\mathcal{T}}(Q)$ is a prime(left) quasi-interior ideal of $\mathcal{M}_1$, if $\overline{\mathcal{T}}(Q) \neq \emptyset$.
(ii) $\underline{\mathcal{T}}(Q)$ is a prime(left) quasi-interior ideal of $\mathcal{M}_1$, if $\underline{\mathcal{T}}(Q) \neq \emptyset$ and $\mathcal{T}$ is strong set-valued anti-homomorphism.

**Proof.** (i) By theorem 5.9.(i), $\overline{\mathcal{T}}(Q)$ is a (left)quasi-interior ideal of $\mathcal{M}_1$. Take aαb $\in \overline{\mathcal{T}}(Q)$ and $\mathcal{T}(a\alpha b) \in Q$. Then $\mathcal{T}(a\alpha b) \cap Q \neq \emptyset$. We have, $\mathcal{T}(a\alpha b) \supseteq \mathcal{T}(b)\alpha\mathcal{T}(a), \forall a, b \in \mathcal{M}_1$ and α $\in \Gamma$. There exists u $\in \mathcal{T}(a), v \in \mathcal{T}(b)$ such that vαu $\subseteq \mathcal{T}(a\alpha b)$. (i.e) vαu $\in Q$. Since $Q$ is prime(left) quasi-interior ideal of $\mathcal{M}_2$, u $\in Q$ or v $\in Q$. Therefore a $\in \overline{\mathcal{T}}(Q)$ or b $\in \overline{\mathcal{T}}(Q)$. Hence $\overline{\mathcal{T}}(Q)$ is a prime(left) quasi-interior ideal of $\mathcal{M}_1$.

(ii) By theorem 5.9.(ii), $\underline{\mathcal{T}}(Q)$ is a (left) quasi-interior ideal of $\mathcal{M}_1$. Suppose $\underline{\mathcal{T}}(Q)$ is not prime(left) quasi-interior ideal of $\mathcal{M}_1$. There exists a, b $\in \mathcal{M}_1$ such that aαb $\in \underline{\mathcal{T}}(Q)$ but a $\notin \underline{\mathcal{T}}(Q)$, b $\notin \underline{\mathcal{T}}(Q)$. Then, as u $\in \mathcal{T}(a), v \in \mathcal{T}(b)$, it follows that u, v $\notin Q$. Thus vαu $\in \mathcal{T}(b)\alpha\mathcal{T}(a) \subseteq \mathcal{T}(a\alpha b) \subseteq Q, \forall a, b \in \mathcal{M}_1$ and α $\in \Gamma$. Since $Q$ is a prime(left) quasi-interior ideal of $\mathcal{M}_2$, u $\in Q$ or v $\in Q$ which is incongruous. Hence $\underline{\mathcal{T}}(Q)$ is a prime(left) quasi-interior ideal of $\mathcal{M}_1$.

**Example 23.** Consider the Γ-semigroup $\mathcal{M}$ of Example 4. Define a set-valued mapping $\mathcal{T}: \mathcal{M} \to \mathcal{P}^*(\mathcal{M})$ by $\mathcal{T}(x_1) = \{x_1, x_2, x_3, x_4\}, \mathcal{T}(x_2) = \{x_1, x_3\}, \mathcal{T}(x_3) = \{x_3\}, \mathcal{T}(x_4) = \{x_4\}$. Then $\mathcal{T}$ is a set-valued anti-homomorphism. $\mathcal{M} = \{x_1, x_2, x_3, x_4\}$ and $\Gamma = \{\alpha\}$ and considering the subsets $\{x_1, x_3\}, \{x_2, x_3\}, \{x_1, x_4\}$. Let $Q = \{x_1, x_2, x_3, x_4\} \subseteq \mathcal{M}$. Then $Q$ is a (left)quasi-interior ideal of $\mathcal{M}$. $\overline{\mathcal{T}}(Q) = \{x_1, x_2, x_3, x_4\}, \mathcal{M}\Gamma Q\Gamma\mathcal{M}\Gamma Q = \{x_1, x_3\}$, $\overline{\mathcal{T}}(\mathcal{M}\Gamma Q\Gamma\mathcal{M}\Gamma Q) = \{x_1, x_2, x_3, x_4\}$. Hence $\overline{\mathcal{T}}(\mathcal{M}\Gamma Q\Gamma\mathcal{M}\Gamma Q) \subseteq \overline{\mathcal{T}}(\mathcal{B})$. Therefore $\overline{\mathcal{T}}(Q)$ is a (left)quasi-interior ideal of $\mathcal{M}$. It is obvious that $Q$ is a prime(left) quasi-interior ideal of Γ-semigroup $\mathcal{M}$. Here, the set $\overline{\mathcal{T}}(Q)$ is a prime(left) quasi-interior ideal of $\mathcal{M}$ for $x_1 \alpha x_4 \in \overline{\mathcal{T}}(Q)$ as $x_1 \in \overline{\mathcal{T}}(Q), x_4 \in \overline{\mathcal{T}}(Q)$.

**Example 24.** Consider the $\Gamma$-semigroups $\mathcal{M}_1$ and $\mathcal{M}_2$ of Example 5. Define a set-valued mapping $\mathcal{T}: \mathcal{M}_1 \to \mathcal{P}^*(\mathcal{M}_2)$ Assume $\mathcal{T}(1) = \{c\}, \mathcal{T}(2) = \{a\}$. Then $\mathcal{T}$ is a strong set-valued anti-homomorphism. Assuming the subsets $\{b\}, \{c\}, \{a, b\}$. Let $Q = \{a, b, c\} \subset \mathcal{M}_2$. Then $Q$ is a quasi-interior ideal of $\mathcal{M}_2$. $\underline{\mathcal{T}}(Q) = \{a, b, c\}$, $\mathcal{M}_2 \Gamma Q \Gamma \mathcal{M}_2 \Gamma Q = \{a, b, c\}$ and $\underline{\mathcal{T}}(\mathcal{M}_2 \Gamma Q \Gamma \mathcal{M}_2 \Gamma Q) = \{a, b, c\}$. Hence $\underline{\mathcal{T}}(\mathcal{M}_2 \Gamma Q \Gamma \mathcal{M}_2 \Gamma Q) \subseteq \underline{\mathcal{T}}(Q)$. Therefore $\underline{\mathcal{T}}(Q)$ is a quasi-interior ideal of $\mathcal{M}_1$. It is clear that $Q$ is a prime quasi-interior ideal of $\Gamma$-semigroup $\mathcal{M}_1$. Also, the set $\underline{\mathcal{T}}(Q)$ is a prime quasi-interior ideal of $\mathcal{M}_1$ for $b\alpha c \in \underline{\mathcal{T}}(Q)$ as $b \in \underline{\mathcal{T}}(Q)$ or $c \in \underline{\mathcal{T}}(Q)$.

**Theorem 5.17.** Consider $\mathcal{M}_1$ and $\mathcal{M}_2$ are two $\Gamma$-semigroups and define $\mathcal{T}: \mathcal{M}_1 \to \mathcal{P}^*(\mathcal{M}_2)$ is a set-valued anti-homomorphism. If $\mathcal{B}$ is a prime bi-quasi-interior ideal of $\mathcal{M}_2$ then

(i)  $\overline{\mathcal{T}}(\mathcal{B})$ is a prime bi-quasi-interior ideal of $\mathcal{M}_1$, if $\overline{\mathcal{T}}(\mathcal{B}) \neq \emptyset$.
(ii) $\underline{\mathcal{T}}(\mathcal{B})$ is a prime bi-quasi-interior ideal of $\mathcal{M}_1$, if $\underline{\mathcal{T}}(\mathcal{B}) \neq \emptyset$ and $\mathcal{T}$ is strong set-valued anti-homomorphism.

**Proof.** (i) By theorem 5.10.(i), $\overline{\mathcal{T}}(\mathcal{B})$ is a bi-quasi-interior ideal of $\mathcal{M}_1$. Assume $a\alpha b \in \overline{\mathcal{T}}(\mathcal{B})$ and $\mathcal{T}(a\alpha b) \in \mathcal{B}$. Then $\mathcal{T}(a\alpha b) \cap \mathcal{B} \neq \emptyset$. We have, $\mathcal{T}(a\alpha b) \supseteq \mathcal{T}(b)\alpha \mathcal{T}(a), \forall a, b \in \mathcal{M}_1$ and $\alpha \in \Gamma$. There exists $u \in \mathcal{T}(a), v \in \mathcal{T}(b)$ such that $v\alpha u \subseteq \mathcal{T}(a\alpha b)$. (i.e) $v\alpha u \in \mathcal{B}$. Since $\mathcal{B}$ is a prime bi-quasi-interior ideal of $\mathcal{M}_2$, $u \in \mathcal{B}$ or $v \in \mathcal{B}$. Thus $a \in \overline{\mathcal{T}}(\mathcal{B})$ or $b \in \overline{\mathcal{T}}(\mathcal{B})$. Hence $\overline{\mathcal{T}}(\mathcal{B})$ is a prime bi-quasi-interior ideal of $\mathcal{M}_1$.

(ii) By theorem 5.10.(ii), $\underline{\mathcal{T}}(\mathcal{B})$ is a bi-quasi-interior ideal of $\mathcal{M}_1$. Suppose $\underline{\mathcal{T}}(\mathcal{B})$ is not prime bi-quasi-interior ideal of $\mathcal{M}_1$. Then there exists $a, b \in \mathcal{M}_1$ such that $a\alpha b \in \underline{\mathcal{T}}(\mathcal{B})$ but $a \notin \underline{\mathcal{T}}(\mathcal{B})$, $b \notin \underline{\mathcal{T}}(\mathcal{B})$. There exists $u \in \mathcal{T}(a), v \in \mathcal{T}(b)$ but $u, v \notin \mathcal{B}$. Thus $v\alpha u \in \mathcal{T}(b)\alpha \mathcal{T}(a) \subseteq \mathcal{T}(a\alpha b) \subseteq \mathcal{B}, \forall a, b \in \mathcal{M}_1$ and $\alpha \in \Gamma$. Since $\mathcal{B}$ is a prime bi-quasi-interior ideal of $\mathcal{M}_2$, $u \in \mathcal{B}$ or $v \in \mathcal{B}$, which is a contradiction. Therefore, $\underline{\mathcal{T}}(\mathcal{B})$ is a prime bi-quasi-interior ideal of $\mathcal{M}_1$.

**Example 25.** Consider the $\Gamma$-semigroup $\mathcal{M}$ of Example 4. Define a set-valued mapping $\mathcal{T}: \mathcal{M} \to \mathcal{P}^*(\mathcal{M})$ by $\mathcal{T}(x_1) = \{x_1, x_2, x_3, x_4\}, \mathcal{T}(x_2) = \{x_1, x_3\}, \mathcal{T}(x_3) = \{x_3\}, \mathcal{T}(x_4) = \{x_4\}$. Then $\mathcal{T}$ is a set-valued anti-homomorphism. $\mathcal{M} = \{x_1, x_2, x_3, x_4\}$ and $\Gamma = \{\alpha\}$ and considering the subsets $\{x_2, x_3\}, \{x_2\}, \{x_1, x_4\}$. Let $\mathcal{B} = \{x_1, x_2, x_3\} \subseteq \mathcal{M}$. Then $\mathcal{B}$ is a bi-quasi-interior ideal of $\mathcal{M}$. $\overline{\mathcal{T}}(\mathcal{B}) = \{x_1, x_2, x_3, x_4\}$, $\mathcal{B}\Gamma\mathcal{M}\Gamma\mathcal{B}\Gamma\mathcal{M}\Gamma\mathcal{B} = \{x_1, x_3\}$, $\overline{\mathcal{T}}(\mathcal{B}\Gamma\mathcal{M}\Gamma\mathcal{B}\Gamma\mathcal{M}\Gamma\mathcal{B}) = \{x_1, x_2, x_3, x_4\}$. Hence $\overline{\mathcal{T}}(\mathcal{B}\Gamma\mathcal{M}\Gamma\mathcal{B}\Gamma\mathcal{M}\Gamma\mathcal{B}) \subseteq \overline{\mathcal{T}}(\mathcal{B})$. Therefore $\overline{\mathcal{T}}(\mathcal{B})$ is a bi-quasi-interior ideal of $\mathcal{M}$. It is obvious that $\mathcal{B}$ is a prime bi-quasi-interior ideal of $\Gamma$-semigroup $\mathcal{M}$. Here, the set $\overline{\mathcal{T}}(\mathcal{B})$ is a prime bi-quasi-interior ideal of $\mathcal{M}$ for $x_1 \alpha x_4 \in \overline{\mathcal{T}}(\mathcal{B})$ as $x_1 \in \overline{\mathcal{T}}(\mathcal{B}), x_4 \in \overline{\mathcal{T}}(\mathcal{B})$.

**Example 26.** Consider the $\Gamma$-semigroups $\mathcal{M}_1$ and $\mathcal{M}_2$ of Example 6. Define a set-valued mapping $\mathcal{T}: \mathcal{M}_1 \to \mathcal{P}^*(\mathcal{M}_2)$ by $\mathcal{T}(x) = \mathcal{T}(z) = \{c\}, \mathcal{T}(y) = \{b, c\}$. Then $\mathcal{T}$ is a strong set-valued anti-homomorphism. Assuming the subsets $\{c\}\{a\}, \{a, b\}$ and Let $\mathcal{B} = \{a, c\} \subset \mathcal{M}_2$. Then $\mathcal{B}$ is a bi-quasi-interior ideal of $\mathcal{M}_2$. Also, $\underline{\mathcal{T}}(\mathcal{B}) = \{a, c\}, \mathcal{B}\Gamma\mathcal{M}\Gamma\mathcal{B}\Gamma\mathcal{M}\Gamma\mathcal{B} = \{a, c\}$ and $\underline{\mathcal{T}}(\mathcal{B}\Gamma\mathcal{M}\Gamma\mathcal{B}\Gamma\mathcal{M}\Gamma\mathcal{B}) = \{a, c\}$. Hence $\underline{\mathcal{T}}(\mathcal{B}\Gamma\mathcal{M}\Gamma\mathcal{B}\Gamma\mathcal{M}\Gamma\mathcal{B}) \subseteq \underline{\mathcal{T}}(\mathcal{B})$. Therefore $\underline{\mathcal{T}}(\mathcal{B})$ is a bi-quasi-interior ideal of $\mathcal{M}_1$. It is obvious that $\mathcal{B}$ is a prime bi-quasi-interior ideal of $\Gamma$- semigroup $\mathcal{M}_1$. Also, the set $\underline{\mathcal{T}}(\mathcal{B})$ is a prime bi-quasi-interior ideal of $\mathcal{M}_1$ for bαc ∈ $\underline{\mathcal{T}}(\mathcal{B})$ as b ∈ $\underline{\mathcal{T}}(\mathcal{B})$ or c ∈ $\underline{\mathcal{T}}(\mathcal{B})$.

### $\mathcal{T}$-Rough Quotient Ideals in Γ-semigroups

Define a set-valued anti-homomorphism $\mathcal{T}: \mathcal{M}_1 \to \mathcal{P}^*(\mathcal{M}_2)$, where $\mathcal{M}_1$ and $\mathcal{M}_2$ are Γ-semigroups.

Suppose $\mathcal{M}_1/\mathcal{T} = \{\mathcal{T}(x)/x \in \mathcal{M}_1\}$. It obviously shows that $\mathcal{M}_1/\mathcal{T}$ is a Γ-semigroup.

**Definition 5.16.** Let $\mathcal{T}: \mathcal{M}_1 \to \mathcal{P}^*(\mathcal{M}_2)$ be a set-valued anti-homomorphism where $\mathcal{M}_1$ and $\mathcal{M}_2$ be Γ-semigroups. Following are $\mathcal{T}$-rough quotient approximations, for $\mathcal{H} \in \mathcal{P}^*(\mathcal{M}_2)$:

$$\underline{\mathcal{T}}(\mathcal{H})/\mathcal{T} = \{\mathcal{T}(x)/\mathcal{T}(x) \subseteq \mathcal{H}\} \text{ and } \overline{\mathcal{T}}(\mathcal{H})/\mathcal{T} = \{\mathcal{T}(x)/\mathcal{T}(x) \cap \mathcal{H} \neq \emptyset\}$$

**Theorem 5.17.** Consider $\mathcal{M}_1$ and $\mathcal{M}_2$ are two Γ-semigroups and define $\mathcal{T}: \mathcal{M}_1 \to \mathcal{P}^*(\mathcal{M}_2)$ a set-valued anti-homomorphism. If $\mathcal{S}$ is a sub-Γ-semigroup of $\mathcal{M}_2$ then

(i)   $\overline{\mathcal{T}}(\mathcal{S})/\mathcal{T}$ is a sub-Γ-semigroup of $\mathcal{M}_1/\mathcal{T}$.

(ii)  $\underline{\mathcal{T}}(\mathcal{S})/\mathcal{T}$ is a sub-Γ-semigroup of $\mathcal{M}_1/\mathcal{T}$, if $\mathcal{T}$ is strong set-valued anti-homomorphism.

**Proof.** (i) Assume $\mathcal{T}(a), \mathcal{T}(b) \in \overline{\mathcal{T}}(\mathcal{S})/\mathcal{T}$. Then $\mathcal{T}(a) \cap \mathcal{S} \neq \emptyset$, $\mathcal{T}(b) \cap \mathcal{S} \neq \emptyset$. There exists $x \in \mathcal{T}(a) \cap \mathcal{S}, y \in \mathcal{T}(b) \cap \mathcal{S}$. By theorem 5.2.(i), $\overline{\mathcal{T}}(\mathcal{S})$ is a sub-Γ-semigroup of $\mathcal{M}_1$. We have xαy ∈ $\mathcal{S}$ and xαy ∈ $\mathcal{T}(aαb)$, ∀ x, y ∈ $\mathcal{M}_1$ and α ∈ Γ. Therefore $\mathcal{T}(aαb) \cap \mathcal{S} \neq \emptyset$ which implies that $\mathcal{T}(aαb) \subseteq \overline{\mathcal{T}}(\mathcal{S})/\mathcal{T}$. Hence $\overline{\mathcal{T}}(\mathcal{S})/\mathcal{T}$ is a sub-Γ-semigroup of $\mathcal{M}_1/\mathcal{T}$.

(ii) Let $\mathcal{T}(a), \mathcal{T}(b) \in \underline{\mathcal{T}}(\mathcal{S})/\mathcal{T}$. Then $\mathcal{T}(a) \subseteq \mathcal{S}, \mathcal{T}(b) \subseteq \mathcal{S}$. By theorem 5.2.(ii), $\underline{\mathcal{T}}(\mathcal{S})$ is a sub-Γ-semigroup of $\mathcal{M}_1$. For all a, b ∈ $\mathcal{M}_1$ and α ∈ Γ, we have $\mathcal{T}(a)α\mathcal{T}(b) \subseteq \mathcal{T}(bαa) = \mathcal{T}(aαb) \subseteq \mathcal{S}\Gamma\mathcal{S} \subseteq \mathcal{S}$. Therefore $\underline{\mathcal{T}}(\mathcal{S})/\mathcal{T}$ is a sub-Γ-semigroup of $\mathcal{M}_1/\mathcal{T}$.

**Theorem 5.18.** Consider two Γ-semigroups $\mathcal{M}_1$ and $\mathcal{M}_2$. Defining $\mathcal{T}: \mathcal{M}_1 \to \mathcal{P}^*(\mathcal{M}_2)$ a set-valued anti-homomorphism. Suppose $\mathcal{J}$ is a left (right, two-sided) ideal of $\mathcal{M}_2$. Then

(i)   $\overline{\mathcal{T}}(\mathcal{J})/\mathcal{T}$ is a left (right, two-sided) ideal of $\mathcal{M}_1/\mathcal{T}$.

(ii)  $\underline{\mathcal{T}}(\mathcal{J})/\mathcal{T}$ is a left (right, two-sided) ideal of $\mathcal{M}_1/\mathcal{T}$, if $\mathcal{T}$ is strong set-valued anti-homomorphism.

**Proof.** (i) Consider $\mathcal{T}(x) \in \overline{\mathcal{T}}(\mathcal{J})/\mathcal{T}$ and $\mathcal{T}(y) \in \mathcal{M}_1/\mathcal{T}$. Then $\mathcal{T}(x) \cap \mathcal{J} \neq \emptyset$, hence x ∈ $\overline{\mathcal{T}}(\mathcal{J})$. By theorem 5.3.(i), $\overline{\mathcal{T}}(\mathcal{J})$ is a left ideal of $\mathcal{M}_1$. We have yαx ∈ $\overline{\mathcal{T}}(\mathcal{J})$, ∀x, y ∈ $\mathcal{M}_1$ and α ∈ Γ. For every m = yαx, then $\mathcal{T}(m) \cap \mathcal{J} \neq \emptyset$. On the other hand, $\mathcal{T}(m) = \mathcal{T}(yαx) = \mathcal{T}(xαy) \supseteq \mathcal{T}(y)α\mathcal{T}(x)$. so, $\mathcal{T}(y)α\mathcal{T}(x) \subseteq \overline{\mathcal{T}}(\mathcal{J})/\mathcal{T}$. Hence $\overline{\mathcal{T}}(\mathcal{J})/\mathcal{T}$ is left ideal of $\mathcal{M}_1/\mathcal{T}$. Similarly, the other cases can be proved.

(ii) Assume $\mathcal{T}(x) \in \underline{\mathcal{T}}(\mathcal{J})/\mathcal{T}$ and $\mathcal{T}(y) \in \mathcal{M}_1/\mathcal{T}$. Then $\mathcal{T}(x) \subseteq \mathcal{J}$ which implies that x ∈ $\underline{\mathcal{T}}(\mathcal{J})$. By theorem 5.3.(ii), $\underline{\mathcal{T}}(\mathcal{J})$ is a left ideal of $\mathcal{M}_1$. We have yαx ∈ $\underline{\mathcal{T}}(\mathcal{J})$, ∀x, y ∈ $\mathcal{M}_1$ and α ∈ Γ. For every m = yαx, we have m ∈ $\underline{\mathcal{T}}(\mathcal{J})$ which means $\mathcal{T}(m) \subseteq \mathcal{J}$. Hence $\mathcal{T}(m) \subseteq \underline{\mathcal{T}}(\mathcal{J})/\mathcal{T}$. On the other hand, from m = yαx, we have $\mathcal{T}(m) = \mathcal{T}(yαx) = \mathcal{T}(xαy) \supseteq$

$\mathcal{T}(y)\alpha\mathcal{T}(x)$. So, $\mathcal{T}(y)\alpha\mathcal{T}(x) \subseteq \underline{\mathcal{T}}(\mathcal{I})/\mathcal{T}$. Therefore $\underline{\mathcal{T}}(\mathcal{I})/\mathcal{T}$ is left ideal of $\mathcal{M}_1/\mathcal{T}$. The process is the same for the other cases.

**Corollary 7.** Define $\mathcal{T}\colon \mathcal{M}_1 \to \mathcal{P}^*(\mathcal{M}_2)$ a set-valued anti-homomorphism where $\mathcal{M}_1$ and $\mathcal{M}_2$ are two $\Gamma$-semigroups. Let $\mathcal{B}$ be a bi-ideal of $\mathcal{M}_2$. Then

(i) $\overline{\mathcal{T}}(\mathcal{B})/\mathcal{T}$ is a bi-ideal of $\mathcal{M}_1/\mathcal{T}$.
(ii) $\underline{\mathcal{T}}(\mathcal{B})/\mathcal{T}$ is a bi-ideal of $\mathcal{M}_1/\mathcal{T}$, if $\mathcal{T}$ is strong set-valued anti-homomorphism.

**Corollary 8.** A mapping $\mathcal{T}\colon \mathcal{M}_1 \to \mathcal{P}^*(\mathcal{M}_2)$ is a set-valued anti-homomorphism where $\mathcal{M}_1$ and $\mathcal{M}_2$ represent $\Gamma$-semigroups. Assume $\mathcal{I}$ is an interior ideal of $\mathcal{M}_2$. Then

(i) $\overline{\mathcal{T}}(\mathcal{I})/\mathcal{T}$ is an interior ideal of $\mathcal{M}_1/\mathcal{T}$.
(ii) $\underline{\mathcal{T}}(\mathcal{I})/\mathcal{T}$ is an interior ideal of $\mathcal{M}_1/\mathcal{T}$, if $\mathcal{T}$ is strong set-valued anti-homomorphism.

**Corollary 9.** Define $\mathcal{T}\colon \mathcal{M}_1 \to \mathcal{P}^*(\mathcal{M}_2)$ a set-valued anti-homomorphism where $\mathcal{M}_1$ and $\mathcal{M}_2$ are $\Gamma$-semigroups. If $\mathcal{Q}$ is a quasi-ideal of $\mathcal{M}_2$ then

(i) $\overline{\mathcal{T}}(\mathcal{Q})/\mathcal{T}$ is a quasi-ideal of $\mathcal{M}_1/\mathcal{T}$.
(ii) $\underline{\mathcal{T}}(\mathcal{Q})/\mathcal{T}$ is a quasi-ideal of $\mathcal{M}_1/\mathcal{T}$, if $\mathcal{T}$ is strong set-valued anti-homomorphism.

**Corollary 10.** Define $\mathcal{T}\colon \mathcal{M}_1 \to \mathcal{P}^*(\mathcal{M}_2)$ a set-valued anti-homomorphism where $\mathcal{M}_1$ and $\mathcal{M}_2$ are two $\Gamma$-semigroups. Let $\mathcal{B}$ be a bi-interior ideal of $\mathcal{M}_2$. Then

(i) $\overline{\mathcal{T}}(\mathcal{B})/\mathcal{T}$ is a bi-interior ideal of $\mathcal{M}_1/\mathcal{T}$.
(ii) $\underline{\mathcal{T}}(\mathcal{B})/\mathcal{T}$ is a bi-interior ideal of $\mathcal{M}_1/\mathcal{T}$, if $\mathcal{T}$ is strong set-valued anti-homomorphism.

**Corollary 11.** A mapping $\mathcal{T}\colon \mathcal{M}_1 \to \mathcal{P}^*(\mathcal{M}_2)$ is a set-valued anti-homomorphism where $\mathcal{M}_1$ and $\mathcal{M}_2$ represent $\Gamma$-semigroups. Assume $\mathcal{B}$ is a bi-quasi ideal of $\mathcal{M}_2$. Then

(i) $\overline{\mathcal{T}}(\mathcal{B})/\mathcal{T}$ is a bi-quasi ideal of $\mathcal{M}_1/\mathcal{T}$.
(ii) $\underline{\mathcal{T}}(\mathcal{B})/\mathcal{T}$ is a bi-quasi ideal of $\mathcal{M}_1/\mathcal{T}$, if $\mathcal{T}$ is strong set-valued anti-homomorphism.

**Corollary 12.** Define $\mathcal{T}\colon \mathcal{M}_1 \to \mathcal{P}^*(\mathcal{M}_2)$ a set-valued anti-homomorphism where $\mathcal{M}_1$ and $\mathcal{M}_2$ are $\Gamma$-semigroups. If $\mathcal{Q}$ is a quasi-interior ideal of $\mathcal{M}_2$ then

(i) $\overline{\mathcal{T}}(\mathcal{Q})/\mathcal{T}$ is a quasi-interior ideal of $\mathcal{M}_1/\mathcal{T}$.
(ii) $\underline{\mathcal{T}}(\mathcal{Q})/\mathcal{T}$ is a quasi-interior ideal of $\mathcal{M}_1/\mathcal{T}$, if $\mathcal{T}$ is strong set-valued anti-homomorphism.

**Corollary 13.** Define $\mathcal{T}\colon \mathcal{M}_1 \to \mathcal{P}^*(\mathcal{M}_2)$ a set-valued anti-homomorphism where $\mathcal{M}_1$ and $\mathcal{M}_2$ are $\Gamma$-semigroups. If $\mathcal{B}$ is a bi-quasi-interior ideal of $\mathcal{M}_2$ then

(i) $\overline{\mathcal{T}}(\mathcal{B})/\mathcal{T}$ is a bi-quasi-interior ideal of $\mathcal{M}_1/\mathcal{T}$.
(ii) $\underline{\mathcal{T}}(\mathcal{B})/\mathcal{T}$ is a bi-quasi-interior ideal of $\mathcal{M}_1/\mathcal{T}$, if $\mathcal{T}$ is strong set-valued anti-homomorphism.

## 6 Application

The proposed work will be implemented in image processing. Homomorphic filtering is the key technique in digital image processing which is represented by the function

$$I(x, y) = f(x, y) * g(x, y)$$

where $I, f, g$ are Image, Illumination and reflectance. In the gray scale image if the boundary region are imprecise set-valued homomorphism identifies the pattern. Based on this pattern we can characterize the selected part of an image.

## 7     Conclusion and future directions

In the present paper, we substituted a universe set by a Γ-semigroup and introduced the notions of $\mathcal{T}$-rough (prime) ideals and $\mathcal{T}$-rough quotient ideals in Γ-semigroup under (strong)set-valued anti-homomorphism. First, we defined the generalized roughness or $\mathcal{T}$-rough sets in Γ-semigroups under (strong)set-valued anti-homomorphism. We defined the notions of $\mathcal{T}$-rough bi-ideal, $\mathcal{T}$-rough quasi-ideal, $\mathcal{T}$-rough interior ideal, $\mathcal{T}$-rough bi-quasi ideal, $\mathcal{T}$-rough bi-interior ideal, $\mathcal{T}$-rough quasi-interior ideal, $\mathcal{T}$-rough bi-quasi-interior ideal, $\mathcal{T}$-rough prime bi-ideal, $\mathcal{T}$-rough prime quasi-ideal, $\mathcal{T}$-rough prime interior ideal, $\mathcal{T}$-rough prime bi-interior ideal, $\mathcal{T}$-rough prime bi-quasi ideal, $\mathcal{T}$-rough prime quasi-interior ideal, $\mathcal{T}$-rough bi-quasi-interior ideal, and $\mathcal{T}$-rough quotient ideals of Γ-semigroups under (strong)set-valued anti-homomorphism. At last, we demonstrated the properties of generalized upper and lower approximations of various (prime)ideals and quotient ideals of Γ-semigroups through (strong)set-valued anti-homomorphism. We have been provided examples for the approximations of various (prime)ideals Γ-semigroup under set-valued anti-homomorphism. We obtain several features of the (strong)set-valued anti-homomorphism of $\mathcal{T}$-Rough set on Γ-semigroups. In future, we will study and analyze the generalized rough weak ideals and tri-ideals under set-valued anti-homomorphism on Γ-semigroup. We hope that this research may provide a powerful tool in image processing. We believe that $\mathcal{T}$-rough ideals offered here will turn out to be more useful in the theory and the applications of rough sets, fuzzy sets and soft sets.

**Acknowledgement**

The authors are highly thankful to the management of Sri Venkateswara College of Engineering and the referees for their grant support and constructive comments and suggestions for improving the paper.